\documentclass{amsart}
\usepackage{amscd}
\usepackage{amsmath}
\usepackage{amsfonts}
\usepackage{amssymb}
\usepackage{amsthm}
\usepackage{epsfig}
\usepackage{graphicx}
\usepackage{color}
\usepackage{amssymb}
\usepackage{esint}
\usepackage{graphicx}
\usepackage{caption}
\usepackage{mathrsfs}
\usepackage{multicol}
\usepackage{color}
\usepackage[T1]{fontenc}
\usepackage{textcomp}
\usepackage[utf8]{inputenc}
\usepackage{enumitem}
\setlist[enumerate]{noitemsep}
\setlist[itemize]{noitemsep}
\setlist[description]{noitemsep}
\usepackage{indentfirst}
\usepackage[toc,page]{appendix}
\usepackage{footmisc}
\def\R{\mathbb R}

\def\N{\mathbb N}
\def\C{\mathbb C}
\def\loc{{\text{\upshape loc}}}
\def\eps{\varepsilon}
\def\n{\|{\hspace{-0.12em}}|}
\DeclareMathOperator\DivGamma{div_{_\Gamma}}
\DeclareMathOperator\DivGammaPrimo{div_{_{\Gamma'}}}
\DeclareMathOperator{\Real}{Re}
\DeclareMathOperator{\Ima}{Im}

\DeclareMathOperator\Tr{Tr}
\newcommand{\cal}[1]{{\mathcal #1}}
\title[The damped wave equation with acoustic boundary conditions...]{The damped wave equation with acoustic boundary conditions and non--locally reacting surfaces}
\author{Alessio Barbieri}
\address[A.~Barbieri]
       {Dipartimento di Ingegneria e Scienze dell'Informazione e Matematica\\
       Università degli Studi dell'Aquila\\
       Via Vetoio - I-67100 L'Aquila ITALY}
\email{alessio.barbieri@graduate.univaq.it}
\author{Enzo Vitillaro}
\address[E.~Vitillaro]
       {Dipartimento di Matematica e Informatica, Universit\`a di Perugia\\
       Via Vanvitelli,1 06123 Perugia ITALY}
\email{enzo.vitillaro@unipg.it}

\date{\today}
\thanks{The work was realized within the auspices of the INdAM -- GNAMPA Projects
{\em Equazioni alle derivate parziali: Problemi e Modelli} (Prot\_U-UFMBAZ-2020-000761), and it was also supported by {\em Progetto Equazione delle onde con condizioni acustiche,  finanziato  con  il Fondo  Ricerca  di Base, 2019, della Universit\`a degli Studi di Perugia} and by {\em Progetti Equazioni delle onde con condizioni iperboliche ed acustiche al bordo,  finanziati  con  i Fondi  Ricerca  di Base 2017 and 2018, della Universit\`a degli Studi di Perugia}.}
\begin{document}
\allowdisplaybreaks
\theoremstyle{plain}
\newtheorem{theorem}{Theorem}[section]
\newtheorem{cor}[theorem]{Corollary}
\newtheorem{lemma}[theorem]{Lemma}
\newtheorem{prop}[theorem]{Proposition}
\newtheorem{defn}[theorem]{Definition}
\theoremstyle{remark}
\newtheorem{oss}{Remark}[section]
\newcommand{\bs}{\boldsymbol}
\begin{abstract} The aim of the paper is to study  the problem
$$
\begin{cases} u_{tt}+du_t-c^2\Delta u=0 \qquad &\text{in
$\R\times\Omega$,}\\
\mu v_{tt}- \DivGamma (\sigma \nabla_\Gamma v)+\delta v_t+\kappa v+\rho u_t =0\qquad
&\text{on
$\R\times \Gamma_1$,}\\
v_t =\partial_\nu u\qquad
&\text{on
$\R\times \Gamma_1$,}\\
\partial_\nu u=0 &\text{on $\R\times \Gamma_0$,}\\
u(0,x)=u_0(x),\quad u_t(0,x)=u_1(x) &
 \text{in $\Omega$,}\\
v(0,x)=v_0(x),\quad v_t(0,x)=v_1(x) &
 \text{on $\Gamma_1$,}
\end{cases}$$
where $\Omega$ is a open domain of $\R^N$ with uniformly $C^r$ boundary ($N\ge 2$, $r\ge 1$),
$\Gamma=\partial\Omega$, $(\Gamma_0,\Gamma_1)$ is a relatively open partition of $\Gamma$ with $\Gamma_0$ (but not $\Gamma_1$) possibly empty.
Here $\DivGamma$ and $\nabla_\Gamma$ denote the
Riemannian divergence and gradient operators on $\Gamma$, $\nu$ is the outward normal
to $\Omega$, the coefficients $\mu,\sigma,\delta, \kappa, \rho$ are suitably regular functions on $\Gamma_1$ with $\rho,\sigma$ and $\mu$ uniformly positive,
$d$ is a suitably regular function in $\Omega$ and $c$ is a positive constant.

In this paper we first study well-posedness in the natural energy space and give regularity results. Hence we study asymptotic stability for solutions when $\Omega$ is bounded, $\Gamma_1$ is connected,  $r=2$, $\rho$ is constant and $\kappa,\delta,d\ge 0$.
\end{abstract}

\maketitle
\section{Introduction and main results}\label{intro}
We deal with the damped wave equation posed in a suitably regular open domain of
$\R^N$, supplied with an acoustic  boundary condition on a part of the boundary and an homogeneous Neumann boundary condition on the (possibly empty) remaining part of it. More precisely we consider the
initial- and boundary-value problem
\begin{equation}\label{1.1}
\begin{cases} u_{tt}+du_t-c^2\Delta u=0 \qquad &\text{in
$\R\times\Omega$,}\\
\mu v_{tt}- \DivGamma (\sigma \nabla_\Gamma v)+\delta v_t+\kappa v+\rho u_t =0\qquad
&\text{on
$\R\times \Gamma_1$,}\\
v_t =\partial_\nu u\qquad
&\text{on
$\R\times \Gamma_1$,}\\
\partial_\nu u=0 &\text{on $\R\times \Gamma_0$,}\\
u(0,x)=u_0(x),\quad u_t(0,x)=u_1(x) &
 \text{in $\Omega$,}\\
v(0,x)=v_0(x),\quad v_t(0,x)=v_1(x) &
 \text{on $\Gamma_1$,}
\end{cases}
\end{equation}
where
where $\Omega\subset\R^N$, $N\ge
2$, is an open domain with boundary $\Gamma=\partial\Omega$ uniformly of class $C^r$ in the sense of \cite{Stein1970}, where
the value of $r\in\N\cup\{\infty\}$ will be further specified when needed, so $r=1$ when nothing is said (in most of the paper we shall take $r=1$ or $r=2$).
We assume $\Gamma=\Gamma_0\cup\Gamma_1$, $\overline{\Gamma_0}\cap\overline{\Gamma_1}=\emptyset$,
so $\Gamma_0$ and $\Gamma_1$ are clopen in $\Gamma$, and $\Gamma_1\not=\emptyset$.
All these properties of $\Omega$, $\Gamma_0$ and $\Gamma_1$ will be more formally stated in \S~\ref{Section 2} and we shall refer to them  as assumption (A0).

Moreover $u=u(t,x)$, $v=v(t,y)$, $t\in\R$, $x\in\Omega$, $y\in\Gamma_1$,
$\Delta=\Delta_x$ denotes the Laplace operator with respect to the space
variable, while $\DivGamma$ and $\nabla_\Gamma$ respectively denote the
Riemannian divergence and gradient operators on $\Gamma$. By
$\nu$ we denote the outward normal to $\Omega$ and $c$ is a fixed positive constant.

Acoustic boundary conditions as those in problem \eqref{1.1} have been introduced by Beale and Rosencrans, for general domains, in \cite{beale, beale2,BealeRosencrans} to model acoustic wave propagation, motivated by  studies in Theoretical Acoustics, see \cite[pp.~259--264]{morseingard}.
In this model  $N=3$, $\Gamma_0=\emptyset$, $\Omega$ is either a bounded or an external domain filled
with a fluid which is at rest but for acoustic wave motion. Since the fluid is assumed to be non-viscous, one denotes by $u$ the
velocity potential, so   $-\nabla u$ is the particle velocity, and $u$ satisfies the wave equation $u_{tt}-c^2\Delta
\phi=0$ in $\R\times\Omega$, where $c>0$ is the sound speed in the fluid.  In the Beale--Rosencrans model, one supposes that $\Gamma$ is not rigid but subject to small oscillations, and that each point of it reacts to the excess pressure of the acoustic wave like a
(possibly) resistive harmonic oscillator or spring, so there is no transverse tension between neighboring points of $\Gamma$, i.e. $\sigma\equiv0$ from the mathematical point of view. These surfaces are called locally reacting in \cite[pp.~259--264]{morseingard}.

After their introduction, acoustic boundary conditions for locally reacting surfaces have been studied in several papers, as for instance \cite{CFL2004}, \cite{FL2006}, \cite{FG2000}, \cite{GGG}, \cite{graber}, \cite{JGPhD}, \cite{JGSB2012}, \cite{KJR2016}, \cite{KT2008}, \cite{LLX2018}, \cite{Maatoug2017}, \cite{mugnolo}. When one dismisses the simplifying assumption that neighboring point do not interact (in the terminology of \cite[p.266]{morseingard}) such surfaces are called of {\em extended reaction}.
We shall call those which react like a membrane {\em non-locally reacting} (other types of reactions can be considered). In this case clearly one has  $\sigma>0$ in \eqref{1.1}.

 The simplest case in which $\sigma$  is constant and  the operator $\DivGamma(\sigma \nabla_\Gamma)$ reduces (up to $\sigma$) to the Laplace--Beltrami operator $\Delta_\Gamma$ was briefly considered in \cite[\S 6]{beale} and then studied in \cite{becklin2019global, FMV2011,  FMV2014, VF2013bis,VF2016, VF2017}.   In these papers, the authors assume $\Gamma_0\neq\emptyset$ and the homogeneous boundary condition on it is replaced by the homogeneous Dirichlet boundary condition.

 For the sake of completeness we would like to mention that  several papers in the literature also deal with the wave equation with porous acoustic boundary conditions, where $\mu\equiv0$ in \eqref{1.1}. For example, we refer to \cite{AN2015bis}, \cite{AN2015}, \cite{BB2017}, \cite{BB2018}, \cite{graber2010}.

 Problem \eqref{1.1}, exactly with the assumptions of the present one but with $d\equiv 0$ was recently studied in \cite{mugnvit}. The aim of the present paper is to focus on the case $d\not\equiv 0$.

The coefficients $\mu$, $\sigma$, $\delta$, $\kappa$, $\rho$ and  $d$ are given real functions, respectively on $\Gamma_1$ and on $\Omega$. Throughout the paper we shall assume that they satisfy the following assumptions, depending on the value of $r$:
\begin{itemize}
	\item[(A1)] $\mu,\,\sigma\in W^{r-1,\infty}(\Gamma_1)$ with $\mu_0:=\text{essinf}_{\Gamma_1}{\mu}>0$ and $\sigma_0:=\text{essinf}_{\Gamma_1}\sigma>0$;
	\item[(A2)] $\rho\in W^{r,\infty}(\Gamma_1)$ with $\rho_0:=\text{essinf}_{\Gamma_1}\rho>0$;
	\item[(A3)] $\delta,\,\kappa\in W^{r-1,\infty}(\Gamma_1)$, $d\in W^{r-1,\infty}(\Omega)$.
\end{itemize}
We notice that in (A1--3) we denoted $W^{\infty,\infty}(\Gamma_1)=\bigcap\limits_{n=1}^\infty W^{n,\infty}(\Gamma_1)=C_b^\infty(\Gamma_1)$.
\begin{footnote}{Here and in the sequel the subscript ``$b$'' in spaces of type $C^r$ means that all derivatives up to order $r$ are (not necessarily uniformly when $r=\infty$) bounded. Clearly  Morrey's Theorem is used.}\end{footnote}

The meaning of the Sobolev spaces used above is the standard one when $\Gamma_1$ is compact, while it will be made precise in  \S~\ref{Section 2} in the non-compact case.

The first aim of the paper is to extend the well--posedness and regularity theory of \cite{mugnvit} to our more general problem.
To state the first of our main results we introduce the phase space
\begin{equation}\label{hcors}
	\mathcal{H}=H^1(\Omega)\times H^1(\Gamma_1)\times L^2(\Omega)\times L^2(\Gamma_1),
\end{equation}
making the reader aware that all functions spaces in the paper are complex, as they are commonly used in Acoustics. On the other hand the corresponding spaces of real-valued functions are trivially invariant under the semigroup associated to \eqref{1.1}.

Our first main result establishes well--posedness for problem \eqref{1.1}.

\begin{theorem}[\bf Well--posedness]\label{wellpo}
	Let assumptions (A0-3) hold. Then, for any choice of data $U_0=(u_0,v_0,u_1,v_1)\in\mathcal{H}$, problem \eqref{1.1} has a unique weak solution\footnote[2]{i.e. a solution defined in a distributional sense, see Definition \ref{weaksoldef}.}
	\begin{equation}
		(u,v)\in C(\mathbb{R};H^1(\Omega)\times H^1(\Gamma_1))\cap C^1(\mathbb{R}; L^2(\Omega)\times L^2(\Gamma_1))
	\end{equation}
	continuously depending on data. Moreover, when we also have\footnote[3]{Here and in the sequel $\Delta$, $\partial_\nu$ and $\DivGamma$ here are taken in suitable distributional sense, see  \S1.5.}
	\begin{equation}
		\begin{split}\label{conditions}
			& u_1\in H^1(\Omega),\quad v_1\in H^1(\Gamma_1),\quad\Delta u_0\in L^2(\Omega),\quad\partial_{\nu}{u_0}_{\textbar\Gamma_1}=v_1,\\ & \partial_{\nu}{u_0}_{\textbar\Gamma_0}=0\quad\text{and}\quad\DivGamma(\sigma\nabla_{\Gamma}v_0)\in L^2(\Gamma_1),
		\end{split}
	\end{equation}
	then
	\begin{equation}\label{reg1}
		(u,v)\in C^1(\mathbb{R};H^1(\Omega)\times H^1(\Gamma_1))\cap C^2(\mathbb{R}; L^2(\Omega)\times L^2(\Gamma_1))
	\end{equation}
	and $\eqref{1.1}_{1}-\eqref{1.1}_{4}$ hold almost everywhere.

Finally, when $\rho(y)\equiv\rho_0>0$, introducing  the energy functional $\mathcal{E}\in C(\mathcal{H})$   by
\begin{multline}\label{energyfunctional}
\mathcal{E}(u,v,w,z)=\frac{\rho_0}{2}\int_{\Omega}|\nabla u|^2+\frac{\rho_0}{2c^2}\int_{\Omega}|w|^2\\
+\frac{1}{2}\int_{\Gamma_1}\sigma|\nabla_{\Gamma}v|_{\Gamma}^2+\frac{1}{2}\int_{\Gamma_1}\mu|z|^2+
\frac{1}{2}\int_{\Gamma_1}\kappa|v|^2,
\end{multline}
solutions satisfy for all $s,t\in\mathbb{R}$ the energy identity
	\begin{equation}\label{energy}
\mathcal{E}(u(\tau),v(\tau),u_t(\tau),v_t(\tau))\Big|_{\tau=s}^{\tau=t}=-\int_{s}^{t}
\left[\int_{\Gamma_1}\delta|v_t|^2+\frac{\rho_0}{c^2}\int_{\Omega}d|u_t|^2\right].
\end{equation}
\end{theorem}
The proof of Theorem~\ref{wellpo} relies on the combination of the abstract well--posedness result in \cite{mugnvit}, on  standard perturbation theory
for linear semigroups and on a suitable characterization of  weak solutions as generalized (or mild) solutions in the semigroup sense.

Our second main results concerns optimal regularity of solutions.
As usual in hyperbolic problems higher regularity
requires corresponding regularity of $\Gamma$ and data, as well as compatibility conditions.
To state next result we set, for $1\le n\le r$,
\begin{equation}\label{hn}
	\mathcal{H}^n:=H^n(\Omega)\times H^n(\Gamma_1)\times H^{n-1}(\Omega)\times H^{n-1}(\Gamma_1),
\end{equation}
and we introduce, when $r=\infty$, the spaces
\begin{equation}\label{cinf}
	\begin{split}
		&C_{L^2}^{\infty}(\overline{\Omega})=\{u\in C^{\infty}(\Omega): D^nu\in L^2_n(\Omega)\,\,\forall n\in\mathbb{N}_0\},\\
		&C_{L^2}^{\infty}(\Gamma_1)=\{v\in C^{\infty}(\Gamma_1): D_{\Gamma}^nv\in L^2_n(\Gamma_1)\,\,\forall n\in\mathbb{N}_0\},
	\end{split}
\end{equation}
where $D$ and $D_\Gamma$ respectively denote the differential on $\Omega$ and the covariant derivative on $\Gamma$,
 while $L^2_n(\Omega)$ and $L^2_n(\Gamma_1)$ stand for the space of $n-$times covariant tensor fields with square integrable norm.

\begin{oss}
Since they look quite uncommon we recall the remarks made in \cite[Remark~1.2.2]{mugnvit} about the space $L^2_n(\Omega)$. They also apply, {\em mutatis mutandis},  to the second one, $L^2_n(\Gamma_1)$.
By Morrey's Theorem all elements of $C^\infty_{L^2}(\overline{\Omega})$ continuously extend, with all their derivatives, to $\overline{\Omega}$, so motivating the notation we used. By the same reason $C^\infty_{L^2}(\overline{\Omega})=\bigcap_{n\in\N_0} H^n(\Omega)$,
so it is a Fréchet space with respect to the associated family of seminorms. Hence the notation $C^\infty(\R;C^\infty_{L^2}(\overline{\Omega}))$
is meaningful in the sense of the G\^{a}teaux derivative (see \cite[pp.~72--74]{Hamilton}), and trivially $C^\infty(\R;C^\infty_{L^2}(\overline{\Omega}))\subseteq C^\infty(\R\times\overline{\Omega})$.
Applying Morrey's Theorem again we have the continuous (and possibly strict) inclusions  $C^\infty_{L^2}(\overline{\Omega})\subseteq C^\infty_b(\overline{\Omega})=C^\infty_b(\Omega)\subseteq C^\infty(\overline{\Omega})$. Trivially $C^\infty_{L^2}(\overline{\Omega})= C^\infty_b(\overline{\Omega})$ when $\Omega$ has finite measure, and $C^\infty_{L^2}(\overline{\Omega})= C^\infty(\overline{\Omega})$ when $\Omega$ is bounded.
\end{oss}

We can then state our optimal regularity result, when $\Gamma$ is at least of class $C^2$.
\begin{theorem}[\bf Optimal regularity]\label{reg}
	Let (A0-3) hold, $r\geq2$ and $n\in\mathbb{N}$ such that $2\leq n\leq r$.
Then, for all $U_0=(u_0,v_0,u_1,v_1)\in\mathcal{H}$, the corresponding weak solution $(u,v)$ of \eqref{1.1} enjoys the further regularity
	\begin{equation}\label{regult}
		(u,v)\in\bigcap_{i=0}^{n}C^i(\mathbb{R};H^{n-i}(\Omega)\times H^{n-i}(\Gamma_1)),
	\end{equation}
if and only if $U_0=(u_0,v_0,u_1,v_1)\in\mathcal{H}^n$ and the compatibility conditions
	\begin{equation}\label{compcond}
		\begin{cases}			\partial_{\nu}u_k&=0\qquad\text{on}\,\,\Gamma_0,\quad\text{for}\,\,k=0,\dots,n-2\\
			\partial_{\nu}u_0&=v_1\qquad\text{on}\,\,\Gamma_1,\\
\mu\partial_{\nu}u_1&=\DivGamma(\sigma\nabla_{\Gamma}v_0)-\delta\partial_{\nu}u_0-\kappa v_0-\rho u_1\qquad \text{on}\,\,\Gamma_1,\\
\mu\partial_{\nu}u_k&=\DivGamma(\sigma\nabla_{\Gamma}\partial_{\nu}u_{k-2})-
\delta\partial_{\nu}u_{k-1}-\kappa\partial_{\nu}u_{k-2}-\rho u_k\qquad\text{on}\,\,\Gamma_1,\\&\text{for}\,\,k=2,\dots,n-2,\quad\text{when $n\geq 4$,}
		\end{cases}
	\end{equation}
	hold, were $u_k$,  for $k=2,\dots, n-2$, was recursively set by $u_{k}=c^2\Delta u_{k-2}-du_{k-1}$.
Moreover, when $r=\infty$, we have
	\begin{equation}\label{regfin}
		(u,v)\in C^{\infty}(\mathbb{R};C_{L^2}^{\infty}(\overline{\Omega})\times C_{L^2}^{\infty}(\Gamma_1)).
	\end{equation}
if and only if  $u_0,u_1\in C_{L^2}^{\infty}(\overline{\Omega})$, $v_0,v_1\in C_{L^2}^{\infty}(\Gamma_1)$ and \eqref{compcond} holds for all $n\in\mathbb{N}$.
\end{theorem}

When $d\equiv 0$ the compatibility conditions \eqref{compcond}  reduce to \cite[(1.9)]{mugnvit}, they being written in a simpler but recursive  (and not closed) form.

The proof of Theorem~\ref{reg} is based on the extension of the standard (when $\Gamma_1$ is compact) regularity theory for elliptic problems on $\Omega$ and $\Gamma_1$ given in  \cite{mugnvit}, see Theorems~\ref{bs}--\ref{rege} below, and on standard semigroup theory.

The second aim of this paper is to study the behavior of solutions of problem  \eqref{1.1} when $\Omega$ is bounded, i.e.
the case studied in \cite[Section~6]{mugnvit}, and $\Gamma_1$ is connected. We shall then keep all assumptions in it, and we shall require that $d\ge 0$,
i.e. that  $du_t$ is a damping term.
More formally in the sequel we shall assume, in addition to (A0--3), that
\begin{itemize}
	\item[(A4)] $\Omega$ is bounded,  $\Gamma_1$ is connected, $r=2$,  $\rho(x)\equiv\rho_0>0$,\, $\delta,\kappa,d\geq0$.
\end{itemize}

As it is clear from the energy identity \eqref{energy}, one can expect asymptotic stability  when the system is damped, i.e.
when  $d\not\equiv 0$ or $\delta\not\equiv 0$. By the contrary, when the system is undamped, i.e., $d\equiv 0$ and $\delta\equiv 0$,
the energy is conserved and one can expect to get pure oscillatory solutions.

The case $d\equiv 0$  was studied in detail in \cite{mugnvit}. We briefly recall in the sequel some results from it. In this case  problem \eqref{1.1} possess three types of  solutions for which $u(t,x)=u(t)$, namely
\renewcommand{\labelenumi}{{\roman{enumi})}}
\begin{enumerate}
\item  $u(t,x)\equiv u_0\in\C$, $v(t,x)\equiv 0$, arising independently on $\kappa$;
\item $u(t,x)=u_1t$, $v(t,x)=u_1v^*(x)$, where $v^*\in H^2(\Gamma_1)$
is the unique solution of the elliptic equation
\begin{equation}\label{1.9}
-\DivGamma(\sigma \nabla_\Gamma v^*)+\kappa v^*+\rho_0=0\qquad \text{on $\Gamma_1$,}
\end{equation}
and $u_1\in\C$, only arising when $\kappa\not\equiv 0$;
\item $u(t,x)\equiv 0$, $v(t,x)\equiv v_0\in \C$, only arising when $\kappa\equiv 0$.
\end{enumerate}
When also $\delta\equiv 0$ a Fourier--type decomposition was obtained in \cite[Theorem~1.3.3 and Corollary~1.3.6 ]{mugnvit}. We refer the interest reader to it.
When $\delta\not\equiv 0$, see  \cite[Corollary~1.3.5]{mugnvit},
denoting by ${\cal H}^{N-1}(\Gamma_1)$ the $N-1$ dimensional Hausdorff measure of $\Gamma_1$ and by $|\Omega|$ the $N$--dimensional Lebesgue measure of $\Omega$, solutions of \eqref{1.1} have the following asymptotic behavior as $t\to\infty$:
\renewcommand{\labelenumi}{{\alph{enumi})}}
\begin{enumerate}
\item when $\kappa\not\equiv 0$
$$\begin{alignedat}{4}
&\nabla u(t)\to 0\quad &&\text{in $[L^2(\Omega)]^N$}, \quad
&&u_t(t)\to \tfrac{\int_\Omega u_1-c^2\int_{\Gamma_1}v_0}{|\Omega|-c^2\int_{\Gamma_1}v^*}\quad &&\text{in $L^2(\Omega)$,}\\
&v_t(t)\to 0\quad &&\text{in $L^2(\Gamma_1)$}, \quad
&& v(t)\to \tfrac{\int_\Omega u_1-c^2\int_{\Gamma_1}v_0}{|\Omega|-c^2\int_{\Gamma_1}v^*}\,v^*\quad &&\text{in $H^1(\Gamma_1)$};
\end{alignedat}$$
\item when $\kappa\equiv 0$
$$\begin{alignedat}{4}
&\nabla u(t)\to 0\,\, && \text{in $[L^2(\Omega)]^N$}, \quad && \quad  u_t(t)\to 0\,\, &&\text{in $L^2(\Omega)$},\quad \\
& v(t)\to \tfrac{c^2\int_{\Gamma_1}v_0-\int_\Omega u_1}{c^2{\cal H}^{N-1}(\Gamma_1)}\quad
&&\text{in $H^1(\Gamma_1)$}, &&\quad  v_t(t)\to 0\quad &&\text{in $L^2(\Gamma_1)$}.
\end{alignedat}
$$
\end{enumerate}
Consequently the energy goes to zero as $t\to\infty$  for all data in $\mathcal{H}$
only when $\kappa\equiv 0$, while when $\kappa\not\equiv 0$ on has to restrict to data such that
$\int_\Omega u_1-c^2\int_{\Gamma_1}v_0=0$. This phenomenon was explained in detail in \cite{mugnvit}.

Our final main result  shows that when $d\not\equiv 0$ the linear inner damping terms stabilizes problem \eqref{1.1}, so the phenomenon described above does not occur anymore.

More precisely it shows that problem \eqref{1.1} is asymptotically stable, i.e. as $t\to\infty$ all solutions converge in the phase space to the stationary solutions of it,
which trivially are:
\renewcommand{\labelenumi}{{\roman{enumi})}}
\begin{enumerate}
\item  $u(t,x)\equiv u_0\in\C$, $v(t,x)\equiv 0$, when $\kappa\not\equiv 0$,
\item  $u(t,x)\equiv u_0\in\C$, $v(t,x)\equiv v_0\in \C$,  when $\kappa\equiv 0$.
\end{enumerate}
Indeed the following result holds.

\begin{theorem}[\bf Asymptotic stability]\label{stabt}
	Let assumptions (A0-4) hold, $d\not\equiv 0$ and $(u,v)$ denote the weak solution of \eqref{1.1} corresponding to data $(u_0,v_0,u_1,v_1)\in\mathcal{H}$.
Then, as  $t\to\infty$,
\renewcommand{\labelenumi}{{\alph{enumi})}}
\begin{enumerate}
\item when $\kappa\not\equiv 0$
\begin{equation}\label{S1}
\begin{alignedat}{4}
& u(t)\to c_1 \quad &&\text{in $H^1(\Omega)$}, \quad
&&u_t(t)\to 0\quad &&\text{in $L^2(\Omega)$,}\\
&v(t)\to 0\quad &&\text{in $H^1(\Gamma_1)$}, \quad
&& v_t(t)\to 0\quad &&\text{in $L^2(\Gamma_1)$};
\end{alignedat}
\end{equation}
where $c_1=\alpha U_0\in\C$ and  $\alpha\in \cal{H}'$ is given by
\begin{equation}\label{alpha}
\alpha U_0=\left(\int _\Omega u_1+ du_0-c^2 \int_{\Gamma_1}v_0\right)\Big{/}\left(\int_\Omega d\right).
\end{equation}
\item When $\kappa\equiv 0$
\begin{equation}\label{S2}
\begin{alignedat}{4}
& u(t)\to c_2\,\, && \text{in $H^1(\Omega)$}, \quad && \quad  u_t(t)\to 0\,\, &&\text{in $L^2(\Omega)$},\quad \\
& v(t)\to c_3\quad
&&\text{in $H^1(\Gamma_1)$}, &&\quad  v_t(t)\to 0\quad &&\text{in $L^2(\Gamma_1)$}.
\end{alignedat}
\end{equation}
where $c_2=\beta U_0\in\C$ and $c_3=\gamma U_0\in\C$,  where $\beta,\gamma\in \cal{H}'$ are  given by
\begin{equation}
\begin{aligned}
\label{beta}
\beta U_0=& \dfrac{\int_{\Gamma_1}\delta\left(\int_\Omega u_1+du_0-c^2\int_{\Gamma_1}v_0\right)+c^2\cal{H}^{N-1}(\Gamma_1)
\int_{\Gamma_1} \mu v_1+\delta v_0+\rho_0u_0}{\int_\Omega d\int_{\Gamma_1}\delta+\rho c^2[\cal{H}^{N-1}(\Gamma_1)]^2},
\\
\gamma U_0=& \dfrac{-\rho_0 \cal{H}^{N-1}(\Gamma_1)\left(\int_\Omega u_1+du_0-c^2\int_{\Gamma_1}v_0\right)+\int_\Omega d
\int_{\Gamma_1} \mu v_1+\delta v_0+\rho_0u_0}{\int_\Omega d\int_{\Gamma_1}\delta+\rho c^2[\cal{H}^{N-1}(\Gamma_1)]^2}.
\end{aligned}
\end{equation}
\end{enumerate}
In both case we have $\mathcal{E}(u(t),v(t),u_t(t),v_t(t))\to 0$ as $t\to\infty$.
\end{theorem}

The proof of Theorem~\ref{stabt} is completely different from that of \cite[Corollary~1.3.5]{mugnvit}, as it relies on two {\em ad--hoc} splittings of the phase space $\cal{H}$, in the two cases a) and b) in the statement.

The paper is simply organized: in Section~\ref{Section 2} we introduce some notation and preliminaries,
while Sections~\ref{Section 3},  \ref{Section 4} and  \ref{Section 5} are respectively devoted to
prove Theorems~\ref{wellpo}, \ref{reg} and  \ref{stabt}.
\section{Notation and preliminaries}\label{Section 2}
\subsection{Common notation.}\label{prel}
We shall denote $\N=\{1,2,\ldots\}$ and $\N_0:=\N\cup\{0\}$.
For any $l\in\N$, $x\in\R^l$ and $\eps>0$ we shall denote by $B_\eps(x)$ the ball of $\R^l$ centered ad $x$ of radius $\eps$.
For $x=(x_1,\ldots,x_l),y=(y_1,\ldots,y_l)\in\C^l$ we shall denote $xy=\sum_{i=1}^lx_iy_i$ and  by $\overline{x}$ the conjugate of $x$.
We shall use the standard notation for (complex) Lebesgue and Sobolev spaces of real order on any open subset $\vartheta$ of $\R^l$, referring to \cite{adams}  for details. We shall also use the standard notation for the $\C^l$-valued version of them. For simplicity
$$\|\cdot\|_{\tau,\vartheta}=\|\cdot\|_{L^\tau(\vartheta)},\quad \|\cdot\|_{\tau}=\|\cdot\|_{\tau,\Omega},\quad
\|\cdot\|_{s,\tau,\vartheta}=\|\cdot\|_{W^{s,\tau}(\vartheta)},\quad \|\cdot\|_{s,\tau}=\|\cdot\|_{s,\tau,\Omega}.
$$

We shall use the standard notations $C^m(\vartheta)$ for $m\in\N\cup\{\infty\}$, and  the subscripts ``$c$" and ``$b$"  will respectively denote subspaces of compactly supported functions and of functions with bounded derivatives  up to order $m$. Moreover for any multiindex $\alpha=(\alpha_1,\ldots, \alpha_l)\in\N_0^l$ we shall use the standard notation $D^\alpha=\partial_1^{\alpha_1}\ldots \partial_l^{\alpha_l}$ and  $|\alpha|=\alpha_1+\ldots+\alpha_l$.

Given a Banach  space $X$ we shall use standard notations for Bochner--Lebesgue and Sobolev Spaces of $X$-valued functions,
we shall denote by $X'$ its dual and by $\langle\cdot,\cdot\rangle_X$ the duality product.  Given another Banach space $Y$,
by ${\cal L}(X,Y)$ we shall indicate the space of bounded linear operators between them.
\subsection{Assumption (A0).}
For the sake of clearness we make precise the structural assumption (A0) mentioned  in \S \ref{intro}:
\renewcommand{\labelenumi}{(A{\arabic{enumi})}}
\begin{enumerate}
\setcounter{enumi}{-1}
\item $\Omega\subset\R^N$, $N\ge
2$, is an open domain with boundary $\Gamma=\partial\Omega=\Gamma_0\cup\Gamma_1$, $\overline{\Gamma_0}\cap\overline{\Gamma_1}=\emptyset$, $\Gamma_1\not=\emptyset$, uniformly of class $C^r$, $r\in\N\cup\{\infty\}$, in the sense of \cite{Stein1970}.
More precisely, following \cite[pp.~423--424]{LeoniSobolev2} and \cite[\S~2.2]{mugnvit} we assume that there exist $\eps_0>0$, $M_i>0$ for $i\in\N_0$, $i\le r$, $N_0\in\N$ and a  countable locally finite open cover $\{\Omega_n\}_n$ of $\Gamma$   such that
\renewcommand{\labelenumii}{({\roman{enumii})}}
\begin{enumerate}
\item if $y\in \Gamma$ then $B_{\eps_0}(y)\subseteq \Omega_n$ for some $n\in\N$;
\item no point of $\R^N$ is contained in more than $N_0$ of the $\Omega_n$'s;
\item for each $n\in\N$ there exists a rigid motion $T_n:\R^N\to\R^N$  and $f_n\in C^r(\R^{N-1})$ with
  $$\|D^\alpha f_n\|_{\infty,\R^{N-1}}\le M_m\quad\text{for all $n\in\N$ and $|\alpha|\le m$, $m\in\N_0$, $m\le r$,}$$
     such that $\Omega_n\cap \Omega=\Omega_n\cap T_n(V_n)$, where
     $$V_n=\{(y',y_N)\in \R^{N-1}\times \R: y_N>f_n(y')\}.$$
\end{enumerate}
\end{enumerate}
In the present paper we shall not directly use assumption (A0). Instead we shall use several results in \cite{mugnvit} which hold only in this setting.
In particular, as shown in \S~2.2 of the quoted paper, any relatively clopen subset $\Gamma'$ of $\Gamma$ is a regularly embedded complete submanifold of $\R^N$.
In  the sequel of this preliminary section $\Gamma'$ will denote such a subset of $\Gamma$ and this geometrical structure will be used.

 \subsection{Geometric notation.} We are now going to recall some notation of geometric nature on $\Gamma$, and hence in $\Gamma'$,
 introduced in \cite[\S~2.3]{mugnvit} (referring to it for more details) when $r$ is possible finite and well--known when $r=\infty$
 (see \cite{Boothby, hebey, jost, taylor}).
We shall respectively denote by $T(\Gamma')$, $T^*(\Gamma')$ and $T^p_q(\Gamma')$ the following  bundles on $\Gamma'$:
the  tangent one, the cotangent one and,  for all $p,q\in\N_0$, the tensor bundle of $p$-times contravariant and $q$-times covariant, or of type $(p,q)$, tensors.
All of them are standardly fiber--wise complexified (see \cite{roman})),  their fibers at $y\in\Gamma'$ being respectively denoted by
$T_y(\Gamma')$, $T_y^*(\Gamma')$ and $T^p_q(T_y(\Gamma'))$. Moreover $\overline{\phantom{a}}$,  $\Real $ and $\Ima $ will respectively stand for conjugation,
real and imaginary parts. We shall also denote $T^0_0(\Gamma')=\Gamma'\times\C$, and we shall use  the standard contraction conventions
$T^p=T^p_0$ and  $T_q=T^0_q$. We also recall that $T(\Gamma')=T^1(\Gamma')$ and $T^*(\Gamma')=T_1(\Gamma')$.

In the sequel we shall consider tensor fields of type $(p,q)$ on $\Gamma'$, that is  $u:\Gamma'\to T^p_q(\Gamma')$ with $u(y)\in T^p_q(T_y(\Gamma'))$
for all $y\in\Gamma'$, so getting tangent and cotangent fields as particular cases.
We shall conventionally identify, as usual, tensor fields of type $(0,0)$ with complex valued functions.
In any coordinate systems the components of a tensor field $u$ with respect to the standard frame field will be denoted by $u^{i_1,\ldots,i_p}_{j_1,\ldots,j_q}$,
with the standard contraction convention.

 As usual we shall denote by $C^m(\Gamma')$  the space of (complex-valued) functions of class $C^m$ on $\Gamma'$ for $m=0,\ldots,r$, and by
$C^{m,p}_q(\Gamma')$ the space of tensor fields which components in any chart are of class $C^m$.
Due to transformation laws for tensors this latter notion is well-defined  only  for $m=0,\ldots,r_1$, where we denote
\begin{equation}\label{2.1}
r_1:=r_1(r,p,q)=\begin{cases}
\infty,\qquad &\text{if $r=\infty$},\\
r,\qquad      &\text{if $r<\infty, \quad p+q=0$},\\
r-1,\qquad      &\text{if $r<\infty, \quad p+q>0$}.\\
\end{cases}
\end{equation}
Trivially $\Gamma$ (and hence $\Gamma'$) inherits from $\R^N$ a Riemannian metric, uniquely extended to an Hermitian one on $T(\Gamma)$, in the sequel denoted by
$(\cdot,\cdot)_\Gamma$, given in local coordinates by $(u,v)_\Gamma=g_{ij}u^i\overline{v^j}$
\begin{footnote}{here and in the sequel we shall use the Einstein summation convention for repeated indices}\end{footnote}, with the $g_{ij}$'s  of class $C^{r-1}$.
The metric $(\cdot,\cdot)_\Gamma$ induces the conjugate-linear (fiber-wise defined) Riesz isomorphism $\flat:T(\Gamma')\to T^*(\Gamma')$, with its inverse
$\sharp$, given by  $\langle \flat u,v\rangle_{T(\Gamma')}=(v,u)_\Gamma$ for all $u,v\in T(\Gamma)$, where
$\langle \cdot,\cdot\rangle_{T(\Gamma')}$ denotes the fiber-wise defined duality pairing.
One then defines the induced bundle metric on $T^*(\Gamma')$ by the formula $(\alpha,\beta)_\Gamma=\langle \alpha,\sharp\beta\rangle_{T(\Gamma')}$
for all $\alpha,\beta\in T^*(\Gamma')$. Consequently one has
\begin{equation}\label{2.2}
  (\alpha,\beta)_\Gamma=(\sharp \beta, \sharp\alpha)_\Gamma,\qquad \text{for all $\alpha,\beta\in T^*(\Gamma')$}.
\end{equation}
More generally $(\cdot,\cdot)_{\Gamma}$ extends to a bundle metric (still denoted by the same symbol) on $T^p_q(\Gamma')$ for any $p,q\in\N_0$ (see \cite[p.~442]{Amann2013} for details).
In the sequel we shall denote, on any $T^p_q(\Gamma')$, $|\cdot|_{\Gamma}^2=(\cdot,\cdot)_\Gamma$.

The natural volume element associated to $(\cdot,\cdot)_\Gamma$ on $\Gamma$ will be denoted by $\omega$  and it is given, in local coordinates,
by $\sqrt gdy_1\wedge\ldots\wedge dy^{N-1}$, where $g:=\det (g_{ij})$.
The Riemannian measure associated to $\omega$ on the Borel $\sigma$--algebra of $\Gamma$
naturally completes to the $\sigma$--algebra of subsets of $\Gamma$ which are
measurable with respect to the Hausdorff measure ${\cal H}^{N-1}$, and $d\cal{H}^{N-1}=\omega$.
All boundary integrals will be refereed, without explicitly using the notation $d\cal{H}^{N-1}$, to this measure. We shall denote by $L^\tau(\Gamma')$,
 the standard (complex) Lebesgue space with respect to ${\cal H}^{N-1}$.
More generally, we say that a tensor field on $\Gamma'$ is ${\cal H}^{N-1}$ measurable provided its components on any chart are measurable and we shall denote,
for $p,q\in\N_0$ and $1\le \tau\le \infty$ by $L^{\tau,p}_q(\Gamma')$ the space of measurable tensor fields $u$ of type $(p,q)$ on $\Gamma'$ , modulo
${\cal H}^{N-1}$--a.e. equivalence, such that $|u|_\Gamma\in L^\tau(\Gamma')$.
We shall denote $\|\cdot\|_{\tau,\Gamma'}=\||\cdot |_\Gamma\|_{L^\tau(\Gamma')}$.
Moreover the notation  $L_{q,\loc}^{\tau,p}(\Gamma')$ we shall have the standard meaning.

We shall denote by $\nabla_\Gamma: C^m(\Gamma')\to C^{m-1,1}(\Gamma')$, $m=1,\ldots, r$,  the Riemannian gradient operator,
defined  by
\begin{equation}\label{rina}
 \nabla_\Gamma u=\sharp d_\Gamma \overline{u}\qquad\text{for $u\in C^1(\Gamma')$,}
\end{equation}
where $d_\Gamma$ stands for the differential on $\Gamma$ (see \cite{sternberg}).
In local coordinates, denoting $(g^{ij})=(g_{ij})^{-1}$, we have
\begin{equation}\label{2.4}
\nabla_\Gamma u=g^{ij}\partial_ju\partial_i \qquad \text{for $u\in C^1(\Gamma')$.}
\end{equation}
We also point out that, by \eqref{2.2}--\eqref{rina},
\begin{equation}\label{2.5}
(\nabla_\Gamma u,\nabla_\Gamma v)_\Gamma=(d_\Gamma u,d_\Gamma v)_\Gamma\qquad \text{for $u,v\in C^1(\Gamma')$.}
\end{equation}
Moreover, when $r\ge 2$, we shall denote by   $\DivGamma:C^{m,1}(\Gamma')\to C^{m-1}(\Gamma')$, $m=1,\ldots,r-1$, the Riemannian divergence operator,
defined for $u\in C^{1,1}(\Gamma')$ by the formula
\begin{equation}\label{2.6}
(\DivGamma u)\omega=d_{\rm ext}(\omega\lrcorner  u),
\end{equation} where $d_{\rm ext}$ and $\lrcorner$ respectively denote the exterior derivative on forms  and the interior product (see \cite[p. 66]{taylor}).
In local coordinates we have $\DivGamma u=g^{-1/2}\partial_i (g^{1/2}u^i)$ for all $u\in C^{1,1}(\Gamma')$.
Hence, integrating by parts on coordinate neighborhoods and using a $C^2$  partition of the unity (see \cite[Theorem~4.1, p.~57]{sternberg}),
for any for $u\in C^1(\Gamma')$, $v\in C^{1,1}(\Gamma')$ such that $uv$ is compactly supported we get
\begin{equation}\label{ridi}
\int_{\Gamma'} (\nabla_\Gamma u,v)_\Gamma=-\int_{\Gamma'}u \DivGamma \overline{v}.
\end{equation}
Finally we shall denote by $ D_\Gamma :C^{m,p}_q(\Gamma')\to C^{m-1,p}_{q+1}(\Gamma')$, $p,q\in \N_0$, $m=1,\ldots,r_1$, the covariant derivative operator
\begin{footnote}{Usually denoted by $\nabla $ in the literature. However, the latter symbol indicates throughout this article the gradient in $\Omega$.}\end{footnote},
referring to \cite{mugnvit} for its definition when $p+q>0$, while $D_\Gamma=d_\Gamma$ on scalar fields.
\subsection{Sobolev spaces on $\Gamma$} Sobolev spaces of real order on $\Gamma$  are classical objects in several cases: when $\Gamma$ is compact, see for example \cite{lionsmagenesIII} and \cite{grisvard};
 when $\Gamma$ is  smooth  and the order is an integer, see  \cite{aubinmanifolds, hebey, jost});
 when $\Gamma$ is  smooth and $\Gamma$ has bounded geometry, see \cite{triebel} and related papers.
All these definitions were unified in \cite[Chapter~3]{mugnvit}, essentially extending the approach in \cite{Amann2013} to the present setting.
In this subsection, which can hence be skipped by the reader only interested in cases listed above, we shall recall some definitions and results from \cite{mugnvit}.

We denote, for $\Gamma'\subseteq \Gamma$, $p,q\in\N_0$, $1\le \tau<\infty$ and $m\in\N_0$, $m\le r_1$, the space $W^{m,\tau,p}_q(\Gamma')$ as  the completion
of $C^{m,\tau,p}_q(\Gamma')=\{u\in C^{m,p}_q(\Gamma'): D^i_\Gamma u\in L^{\tau,p}_{q+i}(\Gamma')\quad\text{for $i=0,\ldots,m$}\}$
with respect to the norm
\begin{equation}\label{2.8}
 \|u\|_{m,\tau,\Gamma'}:=\left( \sum_{i=0}^m\|D_\Gamma ^i u\|_{\tau,\Gamma'}^\tau \right)^{1/\tau}.
\end{equation}
According to standard contraction conventions on tensor orders, $W^{m,\tau}(\Gamma')$ denotes, as usual, the function space.
 We shall also use the standard notation $H^{m,p}_q(\Gamma')=W^{m,2,p}_q(\Gamma')$ and $H^m(\Gamma')=W^{m,2}(\Gamma')$.
The space $W^{m,\tau,p}_q(\Gamma')$ is naturally identified  with the subspace of $L^{\tau,p}_q (\Gamma')$  consisting of the $L^\tau$-limits of
Cauchy sequences in $C^{m,\tau,p}_q(\Gamma')$.
Now for all $u\in W^{m,\tau,p}_q(\Gamma')$, taking $(u_n)_n$ a Cauchy sequence in $X$ that converges to it, one sets $D_{\Gamma'} u$ as the limit in $W^{m-1, \tau,p}_{q+1}(\Gamma')$ of the sequence $D_\Gamma u_n$. Since $D_{\Gamma'} $ is trivially the restriction of $D_\Gamma $, we shall only use the latter notation.
In this way the operator $D_\Gamma$ extends by construction to
\begin{equation}\label{2.9}
D_\Gamma\in {\cal L}(W^{m,\tau,p}_q(\Gamma'); W^{m-1,\tau,p}_{q+1}(\Gamma')),\qquad\text{provided $1\le m\le r_1$,}
\end{equation}
and the norm in $W^{m,\tau,p}_q(\Gamma')$ is still given by \eqref{2.8}.
In particular the norm of $H^1(\Gamma')$
is induced by the inner product
\begin{equation}\label{2.10}
  (u,v)_{1,\Gamma'}=\int_{\Gamma'}u\overline{v}+\int_{\Gamma'}(\nabla_\Gamma u,\nabla_\Gamma v)_\Gamma,
\end{equation}
where $\nabla_\Gamma u$ is  still a.e. given by
\eqref{rina}, so \eqref{2.4}--\eqref{2.5} continue to hold a.e. true.

By standardly setting, for $p$, $q$, $\tau$ and  $m$ as before,  the  space $W^{m,\tau,p}_{q,\loc}(\Gamma')$, the operator $D_\Gamma$ in \eqref{2.9} trivially extends to
$D_\Gamma: W^{m,\tau,p}_{q,\loc}(\Gamma')\to W^{m-1,\tau,p}_{q+1,\loc}(\Gamma')$,
the definition being consistent with possible different values of $m$ and $\tau$.
One then sets, for  $p$, $q$  and $m$ as before,
\begin{equation}\label{2.11}
W^{m,\infty ,p}_q(\Gamma') =\{u\in W^{m,1,p}_{q,\loc}(\Gamma'): \, D_\Gamma^i u \in L^{\infty,p}_{q+i}(\Gamma')\quad \text{for $i=0,\ldots, m$}\},
\end{equation}
 endowed with the norm
\begin{equation}\label{2.12}
\|u\|_{m,\infty,\Gamma'} =\max_{i=0,\ldots,m}\|D_\Gamma^i u\|_{\infty,\Gamma'}.
\end{equation}
The spaces $W^{m,\infty}(\Gamma_1)$ used in assumptions (A1--3) are just the particular case of these spaces when $p=q=0$.

We point out, see \cite[Proposition~3.1.5 and Lemma~3.1.6]{mugnvit}, that the spaces defined above can be simply characterized in local coordinates provided one
uses a suitable atlas for $\Gamma'$ arising from assumption (A0). In the sequel we shall use the following multiplier property,
which is a particular case of \cite[Lemma~3.1.4]{mugnvit} and which extends to the present setting  the analogous property which is well--known when $\Gamma'$ is compact.
\begin{lemma}
For all $p,q\in\N_0$, $1\le \tau\le \infty$,  $m\in\N_0$, $m\le r_1$,
there is $c_1=c_1(p,q,m,\tau)>0$ such that
 $$\|u v\|_{m,\tau,\Gamma'}\le c_1 \|u\|_{m,\infty,\Gamma'} \|v\|_{m,\tau, \Gamma'}
 \quad\text{for all $u\in W^{m,\infty}(\Gamma')$ and  $v\in W^{m,\tau,p}_q(\Gamma')$.}$$
 \label{Multiplierlemma}
\end{lemma}
We now turn to Sobolev spaces of real order. We shall restrict, for simplicity, to the case $\tau=2$ which we shall use in the present paper.
For $s\in (0,r_1)\setminus\N$, $p,q\in\N_0$,  we take $m=[s]$, $\theta=s-m$, so $m<s<m+1$, $0<\theta<1$, and we set
\begin{equation}\label{2.13}
H^{s,p}_q(\Gamma')=\left( H^{m,p}_q(\Gamma'),H^{m+1,p}_q(\Gamma')\right)_{\theta,2},
\end{equation}
where $(\cdot,\cdot)_{\theta,2}$ is the real interpolator functor (se \cite[pp.~39--46]{bergh} or \cite[p.~24]{triebel}).
We point out, see \cite[Lemma~3.1.6]{mugnvit}, that also these spaces  can be simply characterized in the same local coordinates mentioned above.

Compactly supported elements of $C^{r_1,p}_q(\Gamma')$ are dense in $H^{s,p}_q(\Gamma_1)$ for $s\in\R$,  $0\le s\le r_1$, so in particular
$C^r_c(\Gamma')$ is dense in $H^s(\Gamma')$ for $s\in\R$, $0\le s\le r$. This fact, although well--known in the compact setting, may fail in the non--compact one
(see \cite{hebey}) and holds true only thanks to assumption (A0).

Finally we set--up Sobolev spaces of negative order,  in the scalar case, by standardly defining $H^s(\Gamma')=[H^{-s}(\Gamma')]'$ for $s\in\R$, $-r\le s\le 0$,
clearly identifying $L^2(\Gamma')$ with its dual.

\subsection{Operators}\label{subsection2.5}
We now introduce the linear bounded operators to be used in the sequel.
By \cite[Theorem~18.40]{LeoniSobolev2}, the standard trace operator $u\mapsto u_{|\Gamma}$ from $C(\overline{\Omega})$ to $C(\Gamma)$,
when restricted to $C(\overline{\Omega})\cap H^1(\Omega)$, has a unique surjective extension $\Tr\in {\cal L}\left(H^1(\Omega),H^{1/2}(\Gamma)\right)$, with
a bounded right--inverse $R_1$. Moreover, for $m\in\N$, $m\le r$, $\Tr\in {\cal L}\left(H^m(\Omega),H^{m-1/2}(\Gamma)\right)$. We shall denote, as usual, $\Tr u=u_{|\Gamma}$ and, when clear, we shall omit  the subscript $_{|\Gamma}$.
Moreover we shall respectively denote by $u_{|\Gamma_i}$  the restriction of $u_{|\Gamma}$  to  $\Gamma_i$ for $i=0,1$.

Consequently, when $r\geq2$, we can set--up the normal derivative of any
$u\in H^2(\Omega)$ by simply setting $\partial_{\nu}u_{|\Gamma}=(\nabla u\cdot\nu)_{|\Gamma}$, so getting
a linear bounded operator from $H^m(\Omega)$ to $H^{m-3/2}(\Gamma)$.
When $r=1$ for any $u\in H^1(\Omega)$ such that $\Delta u\in L^2(\Omega)$ in the sense of distributions and for any $h\in L^2(\Gamma)$ we say that $\partial_{\nu}u_{|\Gamma}=h$ in  distributional sense if
\begin{equation}\label{lap}
	\int_{\Gamma}hv=\int_{\Omega}\nabla u\nabla v+\int_{\Omega}\Delta u\,v\quad\text{for all}\,\,v\in H^1(\Omega).
\end{equation}
By using the operator  $R_1$ one easily gets that $\partial_{\nu}u_{|\Gamma}$
is unique when it exists. Furthermore, the so defined distributional derivative extends the one defined in the trace sense.
Moreover we shall denote by $\partial_\nu u_{|\Gamma_i}$  the restriction of  $\partial_{\nu}u_{|\Gamma}$ to  $\Gamma_i$, $i=0,1$.

The Riemannian gradient operator $\nabla_\Gamma$ defined in $C^1(\Gamma')$ by \eqref{rina}, when restricted to $C^1(\Gamma')\cap H^1(\Gamma')$, uniquely extends by density to
\begin{equation}\label{truerina}
	\nabla_{\Gamma}\in\mathcal{L}(H^{s+1}(\Gamma');H^{s,1}(\Gamma'))\qquad\text{for $s\in\R$, $0\le s\le r-1$}.
\end{equation}
Also the Riemannian divergence operator $\DivGamma$, defined when $r\ge 2$ in $C^{1,1}(\Gamma)$ by \eqref{2.6},
when restricted to $C^{1,1}(\Gamma')\cap H^{1,1}(\Gamma')$ and hence denoted by $\DivGammaPrimo$, extends by density to
\begin{equation}\label{trueridi33}
\DivGammaPrimo\in\mathcal{L}(H^{s,1}(\Gamma');H^{s-1}(\Gamma'))\quad\text{for $s\in\R$, $1\le s\le r-1$.}
\end{equation}
Moreover, since by density the integration by parts formula \eqref{ridi} holds true, when $r\ge 2$, for all $u\in H^1(\Gamma')$ and $v\in H^{1,1}(\Gamma')$,
taking it as fundamental, we set--up, when $r\ge 1$, $\DivGammaPrimo\in\mathcal{L}(L^{2,1}(\Gamma');H^{-1}(\Gamma'))$ by
\begin{equation}\label{trueridi2}
\langle\text{div}_{\Gamma'}u,v\rangle_{H^1(\Gamma')}=-\int_{\Gamma'}(u,\nabla_{\Gamma}\overline{v})_{\Gamma}
\quad\text{for all $u\in L^{2,1}(\Gamma')$ and $v\in H^{1,1}(\Gamma')$.}
\end{equation}
Clearly the so--defined operator coincides when $r\ge 2$ with the one satisfying \eqref{trueridi33} and consequently, by interpolation, we get that
\begin{equation}\label{trueridi}
\DivGammaPrimo\in\mathcal{L}(H^{s,1}(\Gamma');H^{s-1}(\Gamma'))\quad\text{for $s\in\R$, $0\le s\le r-1$.}
\end{equation}

Mainly to simplify the notation we recall that by identifying  each $u\in H^{s,p}_q(\Gamma_i)$, $i=0,1$, $p,q\in\N_0$, $s\in\R$, $0\le s\le r_1$,
 with its trivial extension to $\Gamma$, we can identify  $H^{s,p}_q(\Gamma_i)$ with its isomorphic image in $H^{s,p}_q(\Gamma)$. As a consequence we have
 the splitting
\begin{equation}\label{split}
	H^{s,p}_q(\Gamma)=H^{s,p}_q(\Gamma_0)\oplus H^{s,p}_q(\Gamma_1),
\end{equation}
which extends by duality to $s\in\R$, $-r_1\le s\le r_1$. We shall consistently make this identification throughout the paper.

Consequently, the operator $\DivGamma$ defined in \eqref{trueridi2} splits to $\DivGamma=(\text{div}_{\Gamma_0},\text{div}_{\Gamma_1})$.
Since $\text{div}_{\Gamma}=\text{div}_{\Gamma_1}$ in $H^{s,1}(\Gamma_1)$,  in the sequel we shall use the simpler notation
$\text{div}_{\Gamma}$ and \eqref{trueridi} reads as
\begin{equation}\label{trueridi3}
\text{div}_{\Gamma}\in\mathcal{L}(H^{s,1}(\Gamma_1);H^{s-1}(\Gamma_1))\quad\text{for}\,\,s\in\mathbb{R},\,0\leq s\leq r-1.
\end{equation}

Hence, by assumption (A1), Lemma~\ref{Multiplierlemma} and interpolation we have
\begin{equation}\label{2.19}\text{div}_{\Gamma}(\sigma\nabla_{\Gamma})\in\mathcal{L}(H^{s+1}(\Gamma_1);H^{s-1}(\Gamma_1))
\quad\text{for $s\in\R$, $0\leq s\leq r-1$.}
\end{equation}
Since by \eqref{ridi} and density we have
\begin{equation}\label{ridinaform}
	\int_{\Gamma_1}-\text{div}_{\Gamma}(\sigma\nabla_{\Gamma}u)v=\int_{\Gamma_1}\sigma(\nabla_{\Gamma}u,\nabla_{\Gamma}v)_{\Gamma}=\int_{\Gamma_1}-\text{div}_{\Gamma}(\sigma\nabla_{\Gamma}v)u.
\end{equation}
for all $u,v\in H^1(\Gamma_1)$, we can extend by transposition the operator in \eqref{2.19} to get the operator
\begin{equation}\label{ridinahilb}
	\text{div}_{\Gamma}(\sigma\nabla_{\Gamma})\in\mathcal{L}(H^{s+1}(\Gamma_1);H^{s-1}(\Gamma_1))\quad\text{for}\,\,0\leq s\leq r-1,
\end{equation}
that we shall use in the sequel.

\section{Well--posedness and weak solutions of $\bs{\eqref{1.1}}$}\label{Section 3}
In order to state our first result we endow the space defined in \eqref{hcors} with the inner product given, for  $V_i=(u_i,v_i,w_i,z_i)$, $i=1,2$, by
\begin{equation}\label{innerp}
	\begin{split}
		\left(V_1,V_2\right)_{\mathcal{H}}=&\int_{\Omega}{\nabla u_1\nabla\overline{u_2}}\,+\int_{\Omega}{u_1\overline{u_2}}\,+\int_{\Gamma_1}{\dfrac{\sigma}{\rho}}\left(\nabla_{\Gamma}v_1,\nabla_{\Gamma}v_2\right)_{\Gamma}\\&+\int_{\Gamma_1}\dfrac{v_1\overline{v_2}}{\rho}\,+\dfrac{1}{c^2}\int_{\Omega}w_1\overline{w_2}\,+\int_{\Gamma_1}\dfrac{\mu}{\rho}z_1\overline{z_2}
	\end{split}
\end{equation}
which is well-defined and equivalent to the standard inner product of $\mathcal{H}$ by (A1-2).
We shall denote $\lVert\cdot\rVert_{\mathcal{H}}=(\cdot,\cdot)^{1/2}_{\mathcal{H}}$ the related norm,
which is trivially equivalent to the standard one. To study problem \eqref{1.1} we formally reduce it to a first order problem by setting $w=u_t$ and $z=v_t$,
so to have
\begin{equation}\label{pbwz}
	\begin{cases}
		u_t-w=0 & \text{in}\,\,\mathbb{R}\times\Omega,\\
		v_t-z=0 & \text{on}\,\,\mathbb{R}\times\Gamma_1,\\
		w_t-dw-c^2\Delta u=0 & \text{in}\,\,\mathbb{R}\times\Omega,\\
		\mu z_t-\text{div}_{\Gamma}(\sigma\nabla_{\Gamma}v)+\delta z+\kappa v+\rho w=0 & \text{on}\,\,\mathbb{R}\times\Gamma_1,\\
		z=\partial_{\nu}u & \text{on}\,\,\mathbb{R}\times\Gamma_1,\\
		\partial_{\nu}u=0 & \text{on}\,\,\mathbb{R}\times\Gamma_0,\\
		u(0,x)=u_0(x),\,\,\,w(0,x)=u_1(x) & \text{in}\,\,\Omega,\\
		v(0,x)=v_0(x),\,\,\,z(0,x)=v_1(x) & \text{on}\,\,\Gamma_1.
	\end{cases}
\end{equation}
Then, we define the unbounded operator $A:D(A)\subset\mathcal{H}\to\mathcal{H}$ by
\begin{equation}\label{domain}
	\begin{split}
		D(A)=\{&(u,v,w,z)\in\left[H^1(\Omega)\times H^1(\Gamma_1)\right]^2:\Delta u\in L^2(\Omega),\\&\partial_{\nu}u_{\textbar\Gamma_0}=0,\,\,
\partial_{\nu}u_{\textbar\Gamma_1}=z,\,\,\text{div}_{\Gamma}(\sigma\nabla_{\Gamma}v)\in L^2(\Gamma_1)\}
	\end{split}
\end{equation}
and
\begin{equation}\label{ad}
	A\left(\begin{matrix}
		u\\v\\w\\z
	\end{matrix}\right)
	=\left(\begin{matrix}
		-w\\-z\\dw-c^2\Delta u\\\frac{1}{\mu}\left[-\text{div}_{\Gamma}(\sigma\nabla_{\Gamma}v)+\delta z+\kappa v+\rho w_{\textbar\Gamma_1}\right]
	\end{matrix}\right),
\end{equation}
where $\Delta u$ is taken in the sense of distributions and $\partial_\nu u_{|\Gamma_0}$,  $\partial_\nu u_{|\Gamma_1}$ and $\DivGamma(\sigma\nabla_\Gamma v)$
  are taken in the sense made precise in \S\ref{subsection2.5}.

So, setting $U=(u,v,w,z)$, problem \eqref{pbwz} can be formally written as
\begin{equation}\label{vectpb}
	U'+AU=0\,\,\text{in}\,\,\mathcal{H},\qquad U(0)=U_0:=(u_0,v_0,w_0,z_0).
\end{equation}
We now recall, for the reader's convenience, the classical definitions of strong and generalized (or mild) solution of \eqref{vectpb} (see \cite[pp.~4 and 105]{pazy} or \cite[Chapter II, pp.~145--150]{EngelNagel}) which we shall use in the sequel.
\begin{defn}\label{sol} We say that
	\begin{itemize}
		\item[i)] $U\in C^1(\mathbb{R},\mathcal{H})$ is a \textit{strong solution} of $U'+AU=0$ if $U(t)\in D(A)$ and $U'(t)+AU(t)=0$ for all $t\in\mathbb{R}$;
		\item[ii)] $U\in C(\mathbb{R},\mathcal{H})$ is a \textit{generalized solution} of $U'+AU=0$ if it is the limit of a sequence of strong solutions of it;
		\item[iii)] $U$ is a \textit{strong} or \textit{generalized solution} of \eqref{vectpb} provided it is a solution of $U'+AU=0$
of the same type and satisfies $U(0)=U_0$.
	\end{itemize}
\end{defn}
We now state our first result on problem \eqref{vectpb}.
\begin{theorem}\label{almostwell}
	Under assumptions (A0-3)
	\begin{itemize}
		\item[i)] the operator $-A$ generates a strongly continous group $\{T(t)\}_{t\in\mathbb{R}}$ on $\mathcal{H}$ and,
for any $U_0\in\mathcal{H}$, problem \eqref{vectpb} possesses a unique generalized solution $U$, given by
		$$U(t)=T(t)U_0,$$
		which is also strong provided $U_0\in D(A)$;
		\item[ii)] recursively defining $D(A^n)$, $n\in\mathbb{N}$, by
		$$D(A^{n+1})=\{U\in D(A^n): AU\in D(A^n)\},$$
		then $D(A^n)$ is an Hilbert space when endowed with the inner product given, for all $V,W\in D(A^n)$, by
		\begin{equation}\label{norm}
			\left(V,W\right)_{D(A^n)}=\sum_{i=0}^{n}{\left(A^iV,A^iW\right)_{\mathcal{H}}};
		\end{equation}
		\item[iii)] for all $n\in\mathbb{N}$ the restriction of the operator $-A$ on $D(A^n)$, with domain $D(A^{n+1})$,
generates a strongly continuous group on it. Consequently, denoting $D(A^0)=\mathcal{H}$, the solution $U$ of \eqref{vectpb} enjoys the further regularity
	\begin{equation}\label{regf}
		U\in\bigcap_{i=0}^{n}C^{n-i}(\mathbb{R}; D(A^i))
	\end{equation}
		if and only if  $U_0\in D(A^n)$.
	\end{itemize}
\end{theorem}

\begin{proof}
We first recall the operator used in the undamped case in \cite{mugnvit}, that is
 $\widetilde{A}:D(\widetilde{A})=D(A)\subset\mathcal{H}\to\mathcal{H}$ given  by
\begin{equation}\label{vitop}
	\widetilde{A} \left(\begin{matrix}
		u\\v\\w\\z
	\end{matrix}\right)=\left(\begin{matrix}
	-w\\-z\\-c^2\Delta u\\\frac{1}{\mu}\left[-\text{div}_{\Gamma}(\sigma\nabla_{\Gamma}v)+\delta z+\kappa v+\rho w_{\textbar\Gamma_1}\right]
\end{matrix}\right).
\end{equation}
By \cite[Theorem 4.1.5-i), p. l]{mugnvit} the operator $-\widetilde{A}$ generates a strongly continuous group
$\{\widetilde{T}(t)\}_{t\in\mathbb{R}}$ on $\mathcal{H}$ and trivially 	$A=\widetilde{A}+B$,
where  $B\in \mathcal{L}(\mathcal{H})$ is given by
\begin{equation}\label{operb}
	B \left(\begin{matrix}
		u\\v\\w\\z
	\end{matrix}\right)=\left(\begin{matrix}0\\0\\dw\\0
\end{matrix}\right).
\end{equation}
The result then follows by standard semigroup theory. We give the details in the sequel, for the reader's convenience.

By \cite[Theorem 1.3, p. 158]{EngelNagel}, the operators  $\pm A$ generate two strongly continuous semigroups, so by
\cite[Chapter II, Generation Theorem for Groups, p. 79]{EngelNagel}, the operator $-A$ generates a strongly continuous group $\{T(t)\}_{t\in\mathbb{R}}$ on $\mathcal{H}$.
The proof of i) is then completed by \cite[p.105]{pazy} or \cite[Chapter II, Proposition~6.4~and~Theorem~6.7,~pp.146--150]{EngelNagel}.

To prove ii) we note that, since by \cite[Chapter II, Theorem~1.4~p.~51]{EngelNagel} the operator $A$ is closed,
one trivially gets by induction that $D(A^n)$ is complete with respect to the norm $\|\cdot\|_{D(A^n)}$ induced by \eqref{norm}.

To prove  iii) we note that, since $-A$ generates a strongly continuous semigroup, by \cite[Chapter II, Theorem 1.10, p. 55]{EngelNagel}
there is  $\Lambda>0$ such that $A+\Lambda I$ is bijective. Moreover
 trivially there is $C_n>0$ such that  $\|(A+\Lambda I)^n \cdot\|_{\cal H}\le C_n \|\cdot\|_{D(A^n)}$, an equivalent norm on $D(A^n)$
 is given by  $\|(A+\Lambda I)^n \cdot \|_{\cal H}$, implicitly used in \cite[Chapter II, Definition~5.1~p.~124]{EngelNagel}.
  Hence, by \cite[Chapter II, Proposition~5.2~p.124]{EngelNagel}, we get that $\pm A$ generate two strongly continuous semigroups on $D(A^n)$
  and then, by the already recalled Generation Theorem for Groups, $-A$ generates a group on it.
  To complete the proof of iii) we simply remark that $U_0\in D(A^n)$ if and only if  $U\in C(\R, D(A^i))$ and
  $U^{(n-i)}=(-A)^{n-i}U$ for $i=0,\ldots,n$.
	\end{proof}

Before proving Theorem~\ref{wellpo} we need to precise what we mean by strong, generalized and weak solutions for problem \eqref{1.1}.
\begin{defn}\label{sgensol}
	A \textit{strong} or \textit{generalized solution} of \eqref{1.1} is simply a couple $(u,v)$ such that $u$ and $v$ are respectively the first and the second components of a solution $U=(u,v,w,z)$ of the same type for \eqref{vectpb}.
\end{defn}
Trivially, strong solutions are also generalized ones. By \eqref{ad} and \eqref{vectpb} a strong solution $(u,v)$ belongs to $C^1(\mathbb{R};H^1(\Omega)\times H^1(\Gamma_1))\cap C^2(\mathbb{R};L^2(\Omega)\times L^2(\Gamma_1))$ and the corresponding strong solution of \eqref{vectpb} is $(u,v,u_t,v_t)$. Hence, given any generalized solution $(u,v)$ and denoting with $U=(u,v,w,z)$ the corresponding strong solution of \eqref{vectpb}, taking a sequence $(u^n,v^n, w^n,z^n)=(u^n,v^n,u_t^n,v_t^n)$ of strong solutions converging to it with respect to the topology of $C(\mathbb{R};\mathcal{H})$, for any $\varphi\in C_c^\infty(\mathbb{R})$ we have
$$\int_{-\infty}^{+\infty}u^n_t\varphi=-\int_{-\infty}^{+\infty}u^n\varphi'\quad\text{in}\,\,L^2(\Omega).$$
Consequently, passing to the limit as $n\to\infty$, one has
$$\int_{-\infty}^{+\infty}w\varphi=-\int_{-\infty}^{+\infty}u\varphi'\quad\text{in}\,\,L^2(\Omega),$$
that is $u_t=w$. By the same argument we prove that $v_t=z$, therefore one gets that $U=(u,v,u_t,v_t)$ and
$$(u,v)\in C(\mathbb{R};H^1(\Omega)\times H^1(\Gamma_1))\cap C^1(\mathbb{R};L^2(\Omega)\times L^2(\Gamma_1)).$$
The corresponding generalized solution $U$ of \eqref{vectpb} is again $(u,v,u_t,v_t)$.

Before defining weak solutions of \eqref{1.1}, following the layout of \cite{mugnvit}, we make precise the meaning of weak solutions of the evolution boundary value problem
\begin{equation}\label{evpb}
	\begin{cases}
		u_{tt}+du_t-c^2\Delta u=0 & \text{in}\,\,\mathbb{R}\times\Omega,\\
		\mu v_{tt}-\text{div}_{\Gamma}(\sigma\nabla_{\Gamma}v)+\delta v_t+\kappa v+\rho u_t=0 & \text{on}\,\,\mathbb{R}\times\Gamma_1,\\
		v_t=\partial_{\nu}u & \text{on}\,\,\mathbb{R}\times\Gamma_1,\\
		\partial_{\nu}u=0 & \text{on}\,\,\mathbb{R}\times\Gamma_0.
	\end{cases}
\end{equation}
\begin{defn}\label{distiddef}
	Given
	$$u\in L^1_{\text{loc}}(\mathbb{R};H^1(\Omega))\cap W^{1,1}_{\text{loc}}(\mathbb{R};L^2(\Omega)),\,\,v\in L^1_{\text{loc}}(\mathbb{R};H^1(\Gamma_1))\cap W^{1,1}_{\text{loc}}(\mathbb{R};L^2(\Gamma_1)),$$
	we say that $(u,v)$ is a \textit{weak solution} of \eqref{evpb} if the distributional identities
	\begin{align}
		\int_{-\infty}^{+\infty}\left[-\int_{\Omega}u_t\varphi_t\,+\int_{\Omega}du_t\varphi\,+c^2\int_{\Gamma_1}v_t\varphi\,-c^2\int_{\Omega}\nabla u\nabla\varphi\right]=0\label{soldeb1}\\
		\int_{-\infty}^{+\infty}\int_{\Gamma_1}\left[-\mu v_t\psi_t+\sigma(\nabla_{\Gamma}v,\nabla_{\Gamma}\overline{\psi})_{\Gamma}+\delta v_t\psi+\kappa v\psi-\rho u\psi_t\right]=0\label{soldeb2}
	\end{align}
	hold for all $\varphi\in C_c^{\infty}(\mathbb{R}\times\mathbb{R}^N)$ and $\psi\in C_c^r(\mathbb{R}\times\Gamma_1)$.
\end{defn}
Trivially, each weak solution of \eqref{evpb} possesses a (unique) representative $(u,v)\in C(\mathbb{R}; L^2(\Omega)\times L^2(\Gamma_1))$. For this reason, in the sequel we shall always consider this representative, so that $u$ and $v$ possess a pointwise meaning. As far as $u_t$ and $v_t$ are concerned, finding their representatives is slightly more difficult and for the reader's convenience we shall provide a proof.
\begin{lemma}\label{reprlemma}
	Under assumptions (A0-3), let $(u,v)$ be a weak solution of \eqref{evpb}. Then,
	\begin{itemize}
		\item[i)] $u_t$ has a unique representative in $C(\mathbb{R};[H^1(\Omega)]')$, which satisfies the alternative distributional identity
		\begin{equation}\label{distid1}
			\begin{split}
				\int_{s}^{t}\left[-\int_{\Omega}u_t\varphi_t\,+c^2\left(\int_{\Omega}\nabla u\nabla\varphi\,-\int_{\Gamma_1}v_t\varphi\right)+\int_{\Omega}du_t\varphi\right]\\+\left[\langle u_t,\varphi\rangle_{H^1(\Omega)}\right]_{s}^{t}=0,
			\end{split}
		\end{equation}
		for all $s,t\in\mathbb{R}$ and $\varphi\in C(\mathbb{R};H^1(\Omega))\cap C^1(\mathbb{R};L^2(\Omega))$;
		\item[ii)] if $u_{\textbar\Gamma_1}\in C(\mathbb{R};H^{-1}(\Gamma_1))$, then $v_t$ possesses a unique representative in $C(\mathbb{R};H^{-1}(\Gamma_1))$ as well, which satisfies the alternative distributional identity
		\begin{equation}\label{distid2}
			\begin{split}
				\int_{s}^{t}\int_{\Gamma_1}\left[-\mu v_t\psi_t+\sigma(\nabla_{\Gamma}v,\nabla_{\Gamma}\overline{\psi})_{\Gamma}+\delta v_t\psi+\kappa v\psi-\rho u\psi_t\right]\\+\left[\langle\mu v_t+\rho u_{\textbar\Gamma_1},\psi\rangle_{H^1(\Gamma_1)}\right]_{s}^{t}=0,
			\end{split}
		\end{equation}
		for all $s,t\in\mathbb{R}$ and $\psi\in C(\mathbb{R};H^1(\Gamma_1))\cap C^1(\mathbb{R};L^2(\Gamma_1))$.
	\end{itemize}	
\end{lemma}
\begin{proof}
	Let $(u,v)$ be a weak solution of \eqref{evpb}. Taking test functions $\varphi$ for \eqref{soldeb1} in the separate form
$\varphi(t,x)=\varphi_1(t)\varphi_0(x)$, with $\varphi_1\in C^{\infty}_c(\mathbb{R})$ and $\varphi_0\in C^{\infty}_c(\mathbb{R}^N)$, we get
	\begin{equation*}
		\int_{-\infty}^{+\infty}\left\{\left(-\int_{\Omega}u_t\varphi_0\right)\varphi_1'
+\left[c^2\left(\int_{\Omega}\nabla u\nabla\varphi_0\,-\int_{\Gamma_1}v_t\varphi_0\right)+\int_{\Omega}du_t\varphi_0\right]\varphi_1\right\}=0.
	\end{equation*}
	Hence, for any fixed $\varphi_0\in C_c^\infty(\mathbb{R}^N)$ the function $t\mapsto\int_{\Omega}u_t(t)\varphi_0$ belongs to
$W^{1,1}_{\text{loc}}(\mathbb{R})$ and
	$$\left(\int_{\Omega}u_t\varphi_0\right)'=c^2\left(-\int_{\Omega}\nabla u\nabla\varphi_0
\,+\int_{\Gamma_1}v_t\varphi_0\right)-\int_{\Omega}du_t\varphi_0\quad\text{a.e. in}\,\,\mathbb{R}.$$
	As a result,
	\begin{equation}\label{tfci}
		\int_{\Omega}u_t(t)\varphi_0\,-\int_{\Omega}u_t(s)\varphi_0=\int_{s}^{t}\left[c^2\left(-\int_{\Omega}\nabla u\nabla\varphi_0\,+\int_{\Gamma_1}v_t\varphi_0\right)-\int_{\Omega}du_t\varphi_0\right]
	\end{equation}
	for a.a. $s,t\in\mathbb{R}$.
Since restrictions of functions in $C_c^\infty(\mathbb{R}^N)$ are dense in $H^1(\Omega)$ (see \cite[Theorem 11.35 p. 330]{LeoniSobolev2}),
taking for all $\varphi_0\in H^1(\Omega)$ a sequence $(\varphi_{0n})_n$ in $C_c^\infty(\mathbb{R}^N)$ such that
 ${\varphi_{0n}}_{|\Omega}\to\varphi_0$ in $H^1(\Omega)$, we get that \eqref{tfci} holds true for a.a. $s,t\in\mathbb{R}$ and all $\varphi_0\in H^1(\Omega)$.
 As a result, denoting by $\mathcal{B}\in\mathcal{L}(H^1(\Omega)\times L^2(\Omega)\times L^2(\Gamma_1);[H^1(\Omega)]')$ the operator defined by
	$$\langle\mathcal{B}(\mathsf{u},\mathsf{v},\mathsf{w}),\varphi_0\rangle_{H^1(\Omega)}
=-\int_{\Omega}\nabla\mathsf{u}\nabla\varphi_0\,+\int_{\Gamma_1}\mathsf{v}\varphi_0\,-\int_{\Omega}\dfrac{d}{c^2}\mathsf{w}\varphi_0,$$
	we have that $(u_t)'=c^2\mathcal{B}(u,v_t,u_t)\in L^1_{\text{loc}}(\mathbb{R};[H^1(\Omega)]')$ weakly.
 Consequently $u_t\in W^{1,1}_{\text{loc}}(\mathbb{R};[H^1(\Omega)]')$, so it has a (unique) representative in $C(\mathbb{R};[H^1(\Omega)]')$. By a standard density argument we get that for any $\varphi\in C^1(\mathbb{R};H^1(\Omega))\subset W^{1,1}_{\text{loc}}(\mathbb{R};H^1(\Omega))$ we have $\langle u_t,\varphi\rangle_{H^1(\Omega)}\in W^{1,1}_{\text{loc}}(\mathbb{R})$ and, by Leibniz rule,
	$$\langle u_t,\varphi\rangle'_{H^1(\Omega)}=\langle u_t',\varphi\rangle_{H^1(\Omega)}+\langle u_t,\varphi_t\rangle_{H^1(\Omega)}\quad\text{a.e. in}\,\,\mathbb{R},$$
	from which \eqref{distid1} follows for test functions $\varphi\in C^1(\mathbb{R};H^1(\Omega))$.
 By standard time regularization then \eqref{distid1} holds for test functions $\varphi\in C(\mathbb{R};H^1(\Omega))\cap C^1(\mathbb{R};L^2(\Omega))$,
 completing the proof of i).
	
	To prove ii), we remark that, taking functions $\psi$ for \eqref{soldeb2} in the separate form $\psi(t,x)=\psi_1(t)\psi_0(x)$ with $\psi_1\in C^\infty_c(\mathbb{R})$ and $\psi_0\in C^r_c(\Gamma_1)$ and using the same arguments as before, we get
	\begin{equation}\label{www}
		\begin{split}
			\int_{\Gamma_1}[\mu v_t(t)+\rho u_{|\Gamma_1}(t)]\psi_0\,-\int_{\Gamma_1}[\mu v_t(s)+\rho u_{\Gamma_1}(s)]\psi_0\\=-\int_{s}^{t}\left[\int_{\Gamma_1}\sigma(\nabla_{\Gamma}v,\nabla_{\Gamma}\overline{\psi_0})_\Gamma\,+\int_{\Gamma_1}\delta v_t\psi_0\,+\int_{\Gamma_1}\kappa v\psi_0\right]
		\end{split}
	\end{equation}
	for all $\psi_0\in C_c^r(\Gamma_1)$ and a.a. $s,t\in\mathbb{R}$. Using the density of $C_c^r(\Gamma_1)$ in $H^1(\Gamma_1)$ we get that \eqref{www} holds for all $\psi_0\in H^1(\Gamma_1)$. Consequently $\mu v_t+\rho u_{|\Gamma_1}\in W^{1,1}_{\text{loc}}(\mathbb{R};H^{-1}(\Gamma_1))$ and $(\mu v_t+\rho u_{|\Gamma_1})'=\text{div}_{\Gamma}(\sigma\nabla_{\Gamma}v)-\delta v_t-\kappa v$ weakly in this space. Hence $v_t$ has a (unique) representative such that $\mu v_t+\rho u_{|\Gamma_1}\in C(\mathbb{R};H^{-1}(\Gamma_1))$. Then, as $u_{|\Gamma_1}\in C(\mathbb{R};H^{-1}(\Gamma_1))$ and $\frac{1}{\mu},\rho\in L^{\infty}(\Gamma_1)$, we have $v_t\in C(\mathbb{R};H^{-1}(\Gamma_1))$. Using the same density argument used before for the space $W^{1,1}(a,b;H^{-1}(\Gamma_1))$ we then get	$\langle\mu v_t+\rho u_{|\Gamma_1},\psi\rangle'_{H^1(\Gamma_1)}=\langle(\mu v_t+\rho u_{|\Gamma_1})',\psi\rangle_{H^1(\Gamma_1)}+\langle\mu v_t+\rho u_{|\Gamma_1},\psi_t\rangle_{H^1(\Gamma_1)}$ a.e. in $\mathbb{R}$, from which \eqref{distid2} follows for test functions $\psi\in C^1(\mathbb{R};H^1(\Gamma_1))$. Finally, by standard time regularization then \eqref{distid2} holds for test functions $\psi\in C(\mathbb{R};H^1(\Omega))\cap C^1(\mathbb{R};L^2(\Omega)).$
\end{proof}
By Lemma \ref{reprlemma} the following definition makes sense.
\begin{defn}\label{weaksoldef}
	For any $U_0\in L^2(\Omega)\times L^2(\Gamma_1)\times[H^1(\Omega)]'\times H^{-1}(\Gamma_1)$ we say that $(u,v)$ is a \textit{weak solution} of \eqref{1.1} if it is a weak solution of \eqref{evpb}, \mbox{$u_{\textbar\Gamma_1}\in C(\mathbb{R};H^{-1}(\Gamma_1))$} and its representative $(u,v)\in C(\mathbb{R};L^2(\Omega)\times L^2(\Gamma_1))\cap C^1(\mathbb{R};[H^1(\Omega)]'\times H^{-1}(\Gamma_1))$ satisfies $(u(0),v(0),u_t(0),v_t(0))=U_0$.
\end{defn}
Trivially, strong solutions of \eqref{1.1} are also weak ones. More generally, also generalized solutions are weak solutions, since Definition \ref{weaksoldef} is stable with respect to convergence in $C(\mathbb{R};H^1(\Omega)\times H^1(\Gamma_1))\cap C^1(\mathbb{R};L^2(\Omega)\times L^2(\Gamma_1))$.  We now present the last preliminary result we need to finally prove the well-posedness of \eqref{1.1}.
\begin{lemma}\label{uniq}
	Under assumptions (A0-3) weak solutions of \eqref{1.1} are unique.
\end{lemma}
\begin{proof}
	By linearity, proving uniqueness reduces to prove that $U_0=0$ implies $(u,v)=0$ in $\mathbb{R}$.
Moreover, when $U_0=0$ the couple $(\hat{u},\hat{v})$ given by $\hat{u}(t)=-u(-t)$, $\hat{v}(t)=v(-t)$ is still a weak solution of \eqref{1.1}
with vanishing data provided $\delta$ and $d$ are respectively replaced by $-\delta$ and $-d$. Thus, we shall just prove that $(u,v)=0$ in $[0,\infty)$.
To prove it, by repeating  the argument in \cite[Proof of Lemma~4.2.5]{mugnvit}, we fix $t>0$ and test functions $\phi$ and $\psi$, depending on $t$, given by
	\begin{equation}\label{phipsiEvans}
		\varphi(s)=\begin{cases}\int_{s}^{t}\overline{u}(\tau)d\tau & \text{if}\,\,s\leq t,\\\overline{u}(t)(t-s)&\text{if}\,\,s\geq t,\end{cases}\qquad\psi(s)=\begin{cases}
			\frac{c^2}{\rho}\int_{s}^{t}\overline{v}(\tau)d\tau &\text{if}\,\,s\leq t,\\\frac{c^2}{\rho}\overline{v}(t)(t-s) &\text{if}\,\,s\geq t,
		\end{cases}
	\end{equation}
	so that, since $1/\rho\in W^{1,\infty}(\Gamma_1)$, by using Lemma~\ref{Multiplierlemma}, we have $\varphi\in C(\mathbb{R};H^1(\Omega))\cap C^1(\mathbb{R};L^2(\Omega))$, $\psi\in C(\mathbb{R};H^1(\Gamma_1))\cap C^1(\mathbb{R};L^2(\Gamma_1))$, with $\varphi_t=-\overline{u}$, $\psi_t=-\frac{c^2}{\rho}\overline{v}$ in $[0,t]$, $\varphi(t)=0$ and $\psi(t)=0$.
	We then apply Lemma \ref{reprlemma} as follows. Since $u_0=u_1=0$ and $v_1=0$, when $s=0$ the sum of distributional identities
\eqref{distid1} and \eqref{distid2} reads as
	\begin{equation}\label{distidmod}
		\begin{split}
			&\int_{0}^{t}\left[u_t\bar{u}\,-c^2\left(\int_{\Omega}\nabla\varphi\nabla\overline{\varphi}\,-\int_{\Gamma_1}v_t\varphi\,+\int_{\Gamma_1}\frac{\mu}{\rho}v_t\overline{v}\right)+\int_{\Omega}du_t\varphi\right.\\&\left.-\dfrac{1}{c^2}\int_{\Gamma_1}\sigma(\nabla_{\Gamma}(\rho\overline{\psi}_t),\nabla_{\Gamma}\overline{\psi})_{\Gamma}\,+\int_{\Gamma_1}\delta v_t\psi\,-\frac{1}{c^2}\int_{\Gamma_1}\kappa\rho\psi\overline{\psi}_t\,+c^2\int_{\Gamma_1}u\overline{v}\right]=0.
		\end{split}
	\end{equation}
	Since $v_0=\psi(t)=0$, integrating by parts, applying  Leibnitz formula $\nabla_\Gamma (\rho\overline{\psi}_t)=\overline{\psi}_t \nabla_\Gamma\rho +\rho\nabla_\Gamma\overline{\psi}_t=
-\frac{c^2}\rho v \nabla_\Gamma \rho+\rho\nabla_\Gamma\overline{\psi}_t$ and taking the real part we get
	\begin{equation*}
		\begin{split}
			&\frac{1}{2}\int_{0}^{t}\dfrac{d}{dt}\left[\lVert u\rVert_2^2-c^2\lVert\nabla\varphi\rVert_2^2+
c^2\int_{\Gamma_1}\dfrac{\mu}{\rho}|v|^2\,-\frac{1}{c^2}\int_{\Gamma_1}\sigma\rho|\nabla_{\Gamma}\psi|_\Gamma^2
-\frac{1}{c^2}\int_{\Gamma_1}\kappa\rho|\psi|^2\right]\\&=-\int_{0}^{t}\int_{\Gamma_1}\dfrac{\sigma}{\rho}\Real \,\left[(\nabla_{\Gamma}\rho,\nabla_{\Gamma}\psi)_\Gamma v\right]-c^2\int_{0}^{t}\int_{\Gamma_1}\dfrac{\delta}{\rho}|v|^2\,-\int_{0}^{t}\int_{\Omega}d|u|^2.
		\end{split}
	\end{equation*}
	Consequently, since $u_0=\varphi(t)=0$ and $v_0=\psi(t)=0$,
	\begin{equation}
		\begin{split}
			&\dfrac{1}{2}\lVert u(t)\rVert_2^2+\dfrac{c^2}{2}\lVert\nabla\varphi(0)\rVert_2^2+\dfrac{c^2}{2}\int_{\Gamma_1}\dfrac{\mu}{\rho}|v(t)|^2\,
+\dfrac{1}{2c^2}\int_{\Gamma_1}\sigma\rho|\nabla_{\Gamma}\psi(0)|_\Gamma^2\\
&=-\dfrac{1}{2c^2}\int_{\Gamma_1}\rho\kappa|\psi(0)|^2\,-\int_{0}^{t}\int_{\Gamma_1}\dfrac{\sigma}{\rho}\Real\,\left[(\nabla_{\Gamma}\rho,\nabla_{\Gamma}\overline{\psi})_\Gamma v\right]\\&-c^2\int_{0}^{t}\int_{\Gamma_1}\dfrac{\delta}{\rho}|v|^2\,-\int_{0}^{t}\int_{\Omega}d|u|^2.
		\end{split}
	\end{equation}
	Since test functions change when $t$ changes, in order to let $t$ vary, we just skip the second term in its left-hand side (so getting an inequality) and we express $\psi$ in terms of the unvarying function $\Upsilon\in C(\mathbb{R};H^1(\Gamma_1))\cap C^1(\mathbb{R};L^2(\Gamma_1))$ defined by $\Upsilon(t):=\frac{c^2}{\rho}\int_{0}^{t}v(\tau)d\tau$, so that $\overline{\psi}(\tau)=\Upsilon(t)-\Upsilon(\tau)$ for $\tau\in[0,t]$.
	
	Proceeding in this way, since $\mu\geq\mu_0>0$, we get
	\begin{equation}\label{distidmod3}
		\begin{split}
			&\dfrac{1}{2}\lVert u(t)\rVert_2^2+\dfrac{c^2\mu_0}{2}\int_{\Gamma_1}\dfrac{|v(t)|^2}{\rho}\,
+\dfrac{1}{2c^2}\int_{\Gamma_1}\sigma\rho|\nabla_{\Gamma}\Upsilon(t)|_\Gamma^2\leq-\dfrac{1}{2c^2}\int_{\Gamma_1}\rho\kappa|\Upsilon(t)|^2\\
&-\int_{0}^{t}\int_{\Gamma_1}\dfrac{\sigma}{\rho}\Real\,[(\nabla_{\Gamma}\rho,\nabla_{\Gamma}\Upsilon(t)-\nabla_{\Gamma}\Upsilon(\tau))_\Gamma v(\tau)]d\tau\,
-c^2\int_{0}^{t}\int_{\Gamma_1}\dfrac{\delta}{\rho}|v|^2\\&-\int_{0}^{t}\int_{\Omega}d|u|^2.
		\end{split}
	\end{equation}
	We now estimate the terms in the right-hand side, using assumptions (A1-3). By H\"{o}lder inequality in time it follows that
	\begin{equation}\label{est1}
	-\dfrac{1}{2c^2}\int_{\Gamma_1}\rho\kappa|\Upsilon(t)|^2\leq c^2\lVert\kappa\rVert_{\infty,\Gamma_1}t\int_{0}^{t}\int_{\Gamma_1}\dfrac{|v|^2}{\rho}
	\end{equation}
	and trivially we have
	\begin{align}
		-c^2\int_{0}^{t}\int_{\Gamma_1}\dfrac{\delta}{\rho}|v|^2&\leq c^2\lVert\delta\rVert_{\infty,\Gamma_1}\int_{0}^{t}\int_{\Gamma_1}\dfrac{|v|^2}{\rho},
\label{est2}\\
		-\int_{0}^{t}\int_{\Omega}d|u|^2&\leq\|d\|_\infty\int_{0}^{t}\lVert u(\tau)\rVert_2^2\,d\tau.\label{est3}
	\end{align}
	Moreover, using Cauchy-Schwarz and Young inequalities, one gets that
	\begin{equation}\label{est4}
		\begin{split}
			&-\int_{0}^{t}\int_{\Gamma_1}\dfrac{\sigma}{\rho}\Real\,[(\nabla_{\Gamma}\rho,\nabla_{\Gamma}\Upsilon(t)-\nabla_{\Gamma}\Upsilon(\tau))_{\Gamma}v(\tau)]d\tau\\&\leq\lVert\nabla_{\Gamma}\rho\rVert_{\infty,\Gamma_1}\lVert\sigma\rVert_{\infty,\Gamma_1}\int_{0}^{t}\int_{\Gamma_1}\dfrac{|v|^2}{\rho}\,+\dfrac{\lVert\nabla_{\Gamma}\rho\rVert_{\infty,\Gamma_1}}{2\rho_0^2}t\int_{\Gamma_1}\sigma\rho|\nabla_{\Gamma}\Upsilon(t)|_{\Gamma}^2\\&+\dfrac{\lVert\nabla_{\Gamma}\rho\rVert_{\infty,\Gamma_1}}{2\rho_0^2}\int_{0}^{t}\int_{\Gamma_1}\sigma\rho|\nabla_{\Gamma}\Upsilon|_{\Gamma}^2.
		\end{split}
	\end{equation}
	Plugging the estimates \eqref{est1}, \eqref{est2}, \eqref{est3} and \eqref{est4} in \eqref{distidmod3} we get that
	\begin{multline}\label{distidmod4}
			\dfrac{1}{2}\lVert u(t)\rVert_2^2+\dfrac{c^2\mu_0}{2}\int_{\Gamma_1}\dfrac{|v(t)|^2}{\rho}\,+\dfrac{1}{2c^2}
\int_{\Gamma_1}\sigma\rho|\nabla_{\Gamma}\Upsilon(t)|_{\Gamma}^2\\
\leq\Big[c^2\lVert\kappa\rVert_{\infty,\Gamma_1}t+c^2\lVert\delta\rVert_{\infty,\Gamma_1}
+\lVert\nabla_{\Gamma}\rho\rVert_{\infty,\Gamma_1}\lVert\sigma\rVert_{\infty,\Gamma_1}\Big]\int_{0}^{t}\int_{\Gamma_1}\dfrac{|v|^2}{\rho}\\
+\dfrac{\lVert\nabla_{\Gamma}\rho\rVert_{\infty,\Gamma_1}}{2\rho_0^2}t\int_{\Gamma_1}\sigma\rho|\nabla_{\Gamma}\Upsilon(t)|_{\Gamma}^2\,
+\dfrac{\lVert\nabla_{\Gamma}\rho\rVert_{\infty,\Gamma_1}}{2\rho_0^2}\int_{0}^{t}\int_{\Gamma_1}\sigma\rho|\nabla_{\Gamma}\Upsilon|_{\Gamma}^2\\
+\|d\|_\infty\int_{0}^{t}\lVert u(\tau)\rVert_2^2d\tau\qquad\text{for all $t\geq0$.}
			\end{multline}
	When $\nabla_{\Gamma}\rho=0$, the estimate \eqref{distidmod4} immediately yields, by applying Gronwall inequality \cite[Appendix B, p. 709]{Evans}, that $(u,v)=0$ in $[0,\infty)$. When $\nabla_{\Gamma}\rho\neq0$, we need to fix $t_1=\frac{\rho_0^2}{2c^2\lVert\nabla_{\Gamma}\rho\rVert_{\infty,\Gamma_1}}$,	so that \eqref{distidmod4} implies
	\begin{multline}
					\dfrac{1}{2}\lVert u(t)\rVert_2^2+\dfrac{c^2\mu_0}{2}\int_{\Gamma_1}\dfrac{|v(t)|^2}{\rho}\,+\dfrac{1}{4c^2}
\int_{\Gamma_1}\sigma\rho|\nabla_{\Gamma}\Upsilon(t)|_{\Gamma}^2\\
\leq\|d\|_\infty\int_{0}^{t}\lVert u(\tau)\rVert_2^2d\tau+\dfrac{\lVert\nabla_{\Gamma}\rho\rVert_{\infty,\Gamma_1}}{2\rho_0^2}\int_{0}^{t}
\int_{\Gamma_1}\sigma\rho|\nabla_{\Gamma}\Upsilon|_{\Gamma}^2\\
+\Big[c^2\lVert\kappa\rVert_{\infty,\Gamma_1}t+c^2\lVert\delta\rVert_{\infty,\Gamma_1}
+\lVert\nabla_{\Gamma}\rho\rVert_{\infty,\Gamma_1}\lVert\sigma\rVert_{\infty,\Gamma_1}\Big]\int_{0}^{t}\int_{\Gamma_1}\dfrac{|v|^2}{\rho}
			\end{multline}
	for $t\in[0,t_1]$. By applying Gronwall inequality again we get $(u,v)=0$ in $[0,t_1)$ and the regularity of $(u,v)$ also gives $(u,v,u_t,v_t)=0$ in $[0,t_1)$. Now, by Definitions \ref{distiddef} and \ref{weaksoldef}, for any $t^*\in\mathbb{R}$, the couple $(u^*_t,v^*_t)$ defined by $u_{t^*}(t)=u(t^*+t)$ and $v_{t^*}(t)=v(t^*+t)$ for all $t\in\mathbb{R}$ is still a weak solution of problem \eqref{1.1} corresponding to initial data $(u(t^*),v(t^*),u_t(t^*),v_t(t^*))$. Hence, previous conclusions allow to state, by induction on $n\in\mathbb{N}$, that $(u,v)=0$ in $[0,nt_1]$ for all $n\in\mathbb{N}$, concluding the proof.
\end{proof}
We can finally prove Theorem \ref{wellpo}.
\begin{proof}[Proof of Theorem \ref{wellpo}]
	By Theorem \ref{almostwell}-i) problem \eqref{1.1} has a unique generalized solution, which is also weak as remarked after
Definition \ref{weaksoldef}.
Moreover, by this Lemma, this solution is also unique among weak solutions.
Continuous dependence on data is an immediate consequence of the strong continuity of the group asserted in Theorem \ref{almostwell}-i).
By \eqref{domain}, \eqref{conditions} implies \eqref{reg1} and hence $\eqref{1.1}_1$ holds true in $L^2(\mathbb{R}\times\Omega)$
(and consequently a.e. in $\mathbb{R}\times\Omega$). By the same type of argument, also $\eqref{1.1}_2-\eqref{1.1}_3$ hold true a.e.
in $\mathbb{R}\times\Gamma_1$ and $\eqref{1.1}_4$ does in $\mathbb{R}\times\Gamma_0$, with respect to the product measure on $\mathbb{R}\times\Gamma$.
	
The energy identity holds for strong solutions. Indeed, by multiplying $\eqref{1.1}_1$ for $u_t$, $\eqref{1.1}_2$ for $v_t$ and then integrating by parts,
 a straightforward calculation leads to \eqref{energy}. Finally, by Theorem \ref{almostwell} $D(A)$ is dense in $\mathcal{H}$
 (see \cite[Generation Theorem 3.8, p. 77]{EngelNagel}), then \eqref{energy} also holds for generalized and, by virtue of Lemma \ref{uniq}, weak solutions.
\end{proof}
\section{Regularity when $r\geq2$} \label{Section 4}
Our second main achievement consists in proving optimal regularity for weak solutions of \eqref{1.1}.
Hence our aim is to improve the regularity estimate \eqref{reg1} and, as usual in hyperbolic problems, an increasing regularity of $\Gamma$ is required,
 thus this section is devoted to the case $r\geq2$. Actually, when $r=1$, more regularity on $\Gamma$ would be meaningless.

 We shall preliminarly recall two regularity results from \cite{mugnvit} which will be used the proof of Theorem \ref{reg}.
 The first one  concerns the operator $\text{div}_{\Gamma}(\sigma\nabla_{\Gamma})$ defined in \eqref{ridinahilb} and hence in particular the Laplace-Beltrami operator $\Delta_{\Gamma}$ when $\sigma\equiv1$.
\begin{theorem}[\bf{\cite[Theorem 5.0.1]{mugnvit}}]\label{bs}
	If (A0-3) hold for any $s\in[-r+1,r-1]$, the operator
	$$B_s:=-\text{div}_{\Gamma}(\sigma\nabla_{\Gamma})+I$$
	is an algebraic and topological isomorphism between $H^{s+1}(\Gamma_1)$ and $H^{s-1}(\Gamma_1)$.
\end{theorem}

The second preliminary result extends a classical regularity estimate (see for example \cite{lionsmagenes1}) for elliptic problems with nonhomogeneous Neumann boundary conditions to the case of noncompact boundaries satisfying assumption (A0). In particular, we shall deal with weak solutions of the classical problem
\begin{equation}\label{ellipticgamma}
	\begin{cases}
		-\Delta u+u=f & \text{in}\,\,\Omega,\\\partial_{\nu}u=\gamma & \text{on}\,\,\Gamma,
	\end{cases}
\end{equation}
where $f\in L^2(\Omega)$ and $\gamma\in L^2(\Gamma)$, that is with $u\in H^1(\Omega)$ such that
\begin{equation}\label{weak}
	\int_{\Omega}\nabla u\nabla\varphi\,+\int_{\Omega}u\varphi\,=\int_{\Omega}f\varphi\,+\int_{\Gamma}\gamma\varphi\qquad\text{for all}\,\,\varphi\in H^1(\Omega).
\end{equation}
\begin{theorem}[\bf{\cite[Theorem 5.0.2]{mugnvit}}]\label{rege}
	Let assumption (A0) hold and $f\in H^{s-2}(\Omega)$, $\gamma\in H^{s-3/2}(\Gamma)$ with $s\in\mathbb{R}$, $2\leq s\leq r$. Then the unique weak solution $u$ of \eqref{weak} belongs to $H^s(\Omega)$. Moreover, there is $c_2=c_2(s,\Omega)>0$ such that
	\begin{equation}
		\lVert u\rVert_{s,2}^{2}\leq c_2\left(\lVert f\rVert_{s-2,2}^2+\lVert\gamma\rVert_{s-3/2,2,\Gamma}^2\right)
	\end{equation}
	for all $f\in H^{s-2}(\Omega)$ and $\gamma\in H^{s-3/2}(\Gamma)$.
\end{theorem}
Theorems \ref{bs} and \ref{rege} give us almost all tools we need to prove Theorem \ref{reg}.
However, one last result is necessary. Keeping in mind the definitions \eqref{hn} and \eqref{cinf},
setting  for $2\leq n\leq r$,	\begin{equation}\label{dndef}
		D_{n-1}=\{(u_0,v_0,u_1,v_1)\in\mathcal{H}^n: \eqref{compcond}\,\,\text{holds}\},
	\end{equation}
and recalling the definition of the subspace $D(A^n)$ of $\mathcal{H}$  recursively given in Theorem \ref{almostwell}-ii),
the following result holds.
\begin{lemma}\label{52}
	If assumptions (A0-3) hold and $r\geq2$, then for all $1\leq n<r$ we have $D(A^n)=D_n$.
	Moreover, the norms $\lVert\cdot\rVert_{D(A^n)}$ (defined in \eqref{norm}) and $\lVert\cdot\rVert_{\mathcal{H}^{n+1}}$ are equivalent on it.
\end{lemma}
\begin{proof}
	At first we more explicitly rewrite the sets $D^d_n$. By identifying $L^2(\Gamma)=L^2(\Gamma_0)\oplus L^2(\Gamma_1)$ (see \eqref{split})
and using \eqref{ridinahilb} together with the Trace Theorem we can introduce
 $\mathcal{E}^n\in\mathcal{L}(\mathcal{H}^{n+1};L^2(\Gamma))$, for $1\leq n<r$, setting  for $U=(u,v,w,z)$,
	\begin{equation}\label{opere}
		\begin{cases}
			\mathcal{E}^1(U)=&(\partial_{\nu}u_{|\Gamma_0},\partial_{\nu}u_{|\Gamma_1}-z)\\
			\mathcal{E}^2(U)=&(\partial_{\nu}w_{|\Gamma_0},\mu\partial_{\nu}w_{|\Gamma_1}-\text{div}_{\Gamma}(\sigma\nabla_{\Gamma}v)+\delta\partial_{\nu}u_{|\Gamma_1}+\kappa v+\rho w_{|\Gamma_1})\\
			\mathcal{E}^3(U)=&(c^2\partial_{\nu}\Delta u_{|\Gamma_0}-\partial_{\nu}(dw)_{|\Gamma_0},\mu c^2\partial_{\nu}\Delta u_{|\Gamma_1}
                             -\mu\partial_{\nu}(dw)_{|\Gamma_1}\\
                             &-\text{div}_{\Gamma}(\sigma\nabla_{\Gamma}\partial_{\nu}u)+\delta\partial_{\nu}w_{|\Gamma_1}+\kappa\partial_{\nu}u_{|\Gamma_1}
                             +\rho c^2\Delta u_{|\Gamma_1}-\rho (dw)_{|\Gamma_1})\\
			\mathcal{E}^4(U)=&(c^2\partial_{\nu}\Delta w_{|\Gamma_0}-c^2\partial_{\nu}(d\Delta u)_{|\Gamma_0}
                            +\partial_{\nu}(d^2w)_{|\Gamma_0},\mu c^2\partial_{\nu}\Delta w_{|\Gamma_1}\\
                            &-\mu c^2\partial_{\nu}(d\Delta u)_{|\Gamma_1}\!+\!\mu\partial_{\nu}(d^2w)_{|\Gamma_1}\!
                            -\!\text{div}_{\Gamma}(\sigma\nabla_{\Gamma}\partial_{\nu}w)\!+\!c^2\delta\partial_{\nu}\Delta u_{|\Gamma_1}\\
                            &-\delta\partial_{\nu}(dw)_{|\Gamma_1}
                            +\kappa\partial_{\nu}w+\rho c^2\Delta w_{|\Gamma_1}-\rho c^2d_{|\Gamma_1}\Delta u_{|\Gamma_1}+\rho (d^2w)_{|\Gamma_1})
		\end{cases}
	\end{equation}
	and recursively, for $i\in\mathbb{N}$ and $4\leq i<r-1$,
	\begin{equation}\label{opereind}
		\mathcal{E}^{i+1}(U)=\mathcal{E}^{i-1}(c^2\Delta u-dw,v,c^2\Delta w-c^2d\Delta u+d^2w,z).
	\end{equation}
We also denote $\mathcal{E}^n=(\mathcal{E}^n_{\Gamma_0},\mathcal{E}^n_{\Gamma_1})$.	
By \eqref{compcond} one trivially gets that
	\begin{equation}\label{dnd}
		D_n=\{U\in\mathcal{H}^{n+1}:\mathcal{E}^i(U)=0\quad\text{for}\,\,i=1,\dots,n\},	
	\end{equation}
	so $D_n$ is a closed subspace of $\mathcal{H}^{n+1}$. We shall endow it with the norm inherited from $\mathcal{H}^{n+1}$.

Since by assumptions (A1--3) $\delta,k,\rho,\frac{1}{\mu}\in W^{r-1,\infty}(\Gamma_1)$ and $d\in W^{r-1,\infty}(\Omega)$,
by well--known multiplier properties in Sobolev spaces on $\Omega$ and by Lemma~\ref{Multiplierlemma}  we have
	\begin{equation}\label{adcont}
		A\in\mathcal{L}(D_n;\mathcal{H}^n)\quad\text{for}\,\,1\leq n<r.
	\end{equation}
We now claim that for all $U\in\mathcal{H}^{n+1}\cap D_1$ and for $1\leq n<r$,
	\begin{equation}\label{char}
		\mathcal{E}^n(U)=0\Longleftrightarrow\mathcal{E}^{n-1}A(U)=0.
	\end{equation}
To prove our claim  we first remark that, since $\partial_{\nu}u_{|\Gamma_1}=z$ in $D_1$, one gets by direct verification
that $\mathcal{E}^2=-(\mathcal{E}^1_{\Gamma_0},\mu\mathcal{E}^1_{\Gamma_1})\cdot A$
in $\mathcal{H}^3\cap D_1$, $\mathcal{E}^3=-\mathcal{E}^2\cdot A$ in $\mathcal{H}^4\cap D_1$ and
$\mathcal{E}^4=-\mathcal{E}^3\cdot A$ in $\mathcal{H}^5\cap D_1$. Moreover, by \eqref{opereind}, one gets that
$\mathcal{E}^{i+2}=\mathcal{E}^i\cdot A^2$ in $\mathcal{H}^{i+3}\cap D_1$ for $3\le i<r-2$, from which our claim follows by induction on $n$.

By recursively applying \eqref{char}  and combining it with \eqref{dnd} it follows that
$D_n=\{U\in\mathcal{H}^{n+1}:\mathcal{E}^1A^i(U)=0\,\,\text{for}\,\,i=1,\dots,n-1\}$	and as a result the  recursive formula
	\begin{equation}\label{dndnew}
		D_{n+1}=\{U\in\mathcal{H}^{n+2}\cap D_n:A(U)\in D_n\}
	\end{equation}
holds for $1\le n<r-2$.	

	We now claim  that $D(A^n)=D_n$ for any $1\leq n<r$. We shall prove our claim by induction, starting from $n=1$. By \eqref{domain} we have $D_1\subset D(A)$.
 To prove the reverse inclusion we take $U=(u,v,w,z)\in D(A)$, so by \eqref{domain} and \eqref{ad}, the quadruplet $(u,v,w,z)$ solves the following system with $(h_1,h_2,h_3,h_4)\in\mathcal{H}$
		\begin{equation}\label{system}
		\begin{cases}
			-w=h_1 & \text{in}\,\,\Omega,\\
			-z=h_2 & \text{on}\,\,\Gamma_1,\\
			dw-c^2\Delta u=h_3 & \text{in}\,\,\Omega,\\
			-\text{div}_{\Gamma}(\sigma\nabla_{\Gamma}v)+\delta z+\kappa v+\rho w_{|\Gamma_1}=\mu h_4 & \text{on}\,\,\Gamma_1,\\
			\partial_{\nu}u=z & \text{on}\,\,\Gamma_1,\\
			\partial_{\nu}u=0 & \text{on}\,\,\Gamma_0,
		\end{cases}
	\end{equation}
so in particular the couple $(u,v)$ solves the system
	\begin{equation*}
		\begin{cases}
			-c^2\Delta u=dh_1+h_3 & \text{in}\,\,\Omega,\\
			-\text{div}_{\Gamma}(\sigma\nabla_{\Gamma}v)+\kappa v=\mu h_4+\delta h_2+\rho h_1 & \text{on}\,\,\Gamma_1,\\
			\partial_{\nu}u=h_2 & \text{on}\,\,\Gamma_1,\\
			\partial_{\nu}u=0 & \text{on}\,\,\Gamma_0.
		\end{cases}
	\end{equation*}
Consequently $u$ is a weak solution of problem \eqref{ellipticgamma} with $f=u+\frac{dh_1+h_3}{c^2}$, $\gamma=0$ on $\Gamma_0$ and $\gamma=h_2$ on $\Gamma_1$,
	where trivially $f\in L^2(\Omega)$ and $\gamma\in H^1(\Gamma_1)$. As regards $v$, it satisfies the equation $B_1v=(1-\kappa)v+\mu h_4+\delta h_2+\rho {h_1}_{|\Gamma_1}$,
	where $B_1$ is the operator defined in Theorem \ref{bs}.
Since $f\in L^2(\Omega)$, $\gamma\in H^1(\Gamma_1)$ and trivially $B_1v\in L^2(\Gamma_1)$, by Theorem \ref{bs} we get $v\in H^2(\Gamma_1)$ and
\mbox{applying} Theorem \ref{rege} with $s=2$ we get $u\in H^2(\Omega)$, so that $U\in\mathcal{H}^2$. Finally, by $\eqref{system}_5-\eqref{system}_6$,
$U\in D_1$, so concluding the proof of our claim in the case $n=1$.
	
	We now suppose by induction that $D(A^n)=D_n$ for $n<r-1$. By \eqref{dndnew} we immediately get that
$D_{n+1}\subset D(A^{n+1})$. To prove the reverse inclusion let $U=(u,v,w,z)\in D(A^{n+1})$, so that by the induction hypothesis
$U,A(U)\in D_n\subset\mathcal{H}^{n+1}$. Repeating the same arguments used in the case $n=1$, we get that $u$ and $v$ are solutions of the same problems,
with the only difference that $f\in H^n(\Omega)$, $\gamma\in H^{n+1}(\Gamma_1)$ and $B_1v\in H^n(\Gamma_1)$, thanks to assumptions (A1-3).
 Hence, applying Theorems 3.7 and 3.8 once again, $U\in\mathcal{H}^{n+2}$. Having $U,A(U)\in D_n$ already, $U\in D_{n+1}$, concluding the proof of our last claim.
	
	Now, to prove the equivalence of the norms we stated it is enough to remark that that from \eqref{adcont} it follows that
$A^i\in\mathcal{L}(D(A);\mathcal{H})\quad\text{for all}\,\,i=1,\dots,n$. Combining this fact with \eqref{norm} we get that there exists \mbox{$c_{3}=c_{3}(n,\sigma,\delta,\kappa,\rho,\mu,c,d)>0$} such that
	\begin{equation}\label{equiv3}
		\lVert U\rVert_{D(A^n)}\leq c_{3}\lVert U\rVert_{\mathcal{H}^{n+1}}\quad\text{for all}\,\,U\in D(A^n).
	\end{equation}
	The reverse inequality follows by \cite[Corollary 2.8, p. 35]{brezis2}, since both $\lVert\cdot\rVert_{D(A^n)}$ and
$\lVert\cdot\rVert_{\mathcal{H}^{n+1}}$ are two norms on $D_n$ with respect to which $D_n$ is complete  and \eqref{equiv3} holds.
\end{proof}

Finally, we have all the necessary notions to prove Theorem \ref{reg}.
\begin{proof}[Proof of Theorem \ref{reg}.] The first part of the statement follows by Theorem \ref{almostwell}-iii) and Lemma \ref{52}.
Indeed the solution $U$ of \eqref{vectpb} enjoys the regularity
	\begin{equation}\label{regolu}
		U\in\bigcap_{i=0}^n C^{n-i}(\mathbb{R};D(A^i))
	\end{equation}
if and only if  $U_0\in D(A^n)$, and $D(A^n)=D_n$ by Lemma~\ref{52}.
	The second part immediately follows by the first one: when $r=\infty$ we have no upper bound to $n$ in \eqref{regult}, which thus holds for all $i\in\mathbb{N}$. Furthermore, by Morrey's Theorem
	$$C_{L^2}^\infty(\overline{\Omega})=\bigcap_{n\in\mathbb{N}}H^n(\Omega),\quad C_{L^2}^\infty(\Gamma_1)=\bigcap_{n\in\mathbb{N}}H^n(\Gamma_1),$$
	so if $u_0,u_1\in C_{L^2}^\infty(\overline{\Omega})$ and $v_0,v_1\in C_{L^2}^\infty(\Gamma_1)$, then \eqref{regult} holds.
\end{proof}
\section{Asymptotic stability} \label{Section 5}
This section  is devoted to study the asymptotic stability for problem \eqref{1.1} when assumptions (A0--4) hold, i.e. to prove Theorem~\ref{stabt}.
We start by introducing on $\cal{H}$  the pseudo--inner product $[\cdot,\cdot]_\cal{H}$ suggested by the energy identity \eqref{energy}.
For $U_i=(u_i,v_i,w_i,z_i)\in\cal{H}$, $i=1,2$, we set
\begin{equation}\label{5.1}
		[U_1,U_2]_\cal{H}=\rho_0\!\int_{\Omega}\!\!\!{\nabla u_1\nabla\overline{u_2}}+\!\int_{\Gamma_1}\!\!\!\!\sigma\left(\nabla_{\Gamma}v_1,\nabla_{\Gamma}v_2\right)_{\Gamma}
+\!\!\int_{\Gamma_1}\!\!\!\!\kappa v_1\overline{v_2}+\dfrac{\rho_0}{c^2}\int_{\Omega}\!\!\!\!w_1\overline{w_2}+\!\int_{\Gamma_1}\!\!\!\!\mu z_1\overline{z_2}.
\end{equation}
By assumptions (A0--4) the sesquilinear form $[\cdot,\cdot]_\cal{H}$ is trivially continuous.
We shall respectively denote by
\begin{equation}\label{5.2}
  \n\cdot\n_{\cal{H}}=[\cdot,\cdot]_{\cal{H}}^{1/2},\quad \text{and}\quad\cal{N}=\{U\in \cal{H}: \n U\n_{\cal{H}}=0\},
\end{equation}
the pseudo--norm and the null-- space associated to $[\cdot,\cdot]_\cal{H}$.
We immediately characterize $\cal{N}$.
\begin{lemma}\label{lemma 5.1} Let assumptions (A0--4) hold. Then
\begin{equation}\label{5.3}
\cal{N}=
\begin{cases}
\{(c_1,0,0,0), c_1\in\C\} &\quad\text{when $\kappa\not\equiv 0$},\\
\{(c_1,c_2,0,0), c_1, c_2\in\C\} &\quad\text{when $\kappa\equiv 0$.}
\end{cases}
\end{equation}
\end{lemma}
\begin{proof} Trivially $\cal{N}$ contains the space in the right--hand side of \eqref{5.3}. Moreover, by assumptions (A1--2) and (A4) and by \eqref{5.1},
if $U=(u,v,w,z)$ and $\n U\n=0$ we have   $\nabla u=0$, $\nabla_\Gamma v=0$, $w=0$ and $z=0$.
Since both $\Omega$ and $\Gamma_1$ are connected then $u=c_1\in\C$ in $\Omega$ and $v=c_2\in\C$ on $\Gamma_1$. When $\kappa\not\equiv 0$
we also get $c_2=0$.
\end{proof}
The following key result collects several properties of the operator $A$ in connection with $[\cdot,\cdot]_\cal{H}$.
\begin{lemma}\label{lemma5.2} Let assumptions (A0--4) hold. Then:
\renewcommand{\labelenumi}{{\roman{enumi})}}
\begin{enumerate}
\item $A$ has compact resolvent;
\item $\text{Ker } A=\cal{N}$;
\item for all $U=(u,v,w,z)\in D(A)$ we have
\begin{equation}\label{5.4}
\Real \, [AU,U]_{\cal{H}}=\dfrac {\rho_0}{c^2}\int_\Omega d|w|^2+\int_{\Gamma_1}\delta |z|^2\ge 0;
\end{equation}
\item if $d\not\equiv 0$   for all $\lambda\in\R$ and  $U\in D(A)$ such that $AU=i\lambda U$ we have $U\in\cal{N}$.
\end{enumerate}
\end{lemma}
\begin{proof}We first prove i). Since $r\ge 2$, by Lemma~\ref{52} we have
\begin{equation}\label{5.5}
D(A)=\{(u,v,w,z)\in \cal{H}^2: \partial_\nu u_{|\Gamma_0}=0\quad\text{and}\quad \partial_\nu u_{|\Gamma_1}=z\},
\end{equation}
$\cal{H}^2$ being given by \eqref{hn}, and $\|\cdot\|_{D(A)}$ is equivalent to $\|\cdot\|_{\cal{H}^2}$ on it.
Since $\Omega$ is bounded the injection $\cal{H}^2\hookrightarrow \cal{H}$ is compact, hence by \eqref{5.5} the canonical injection
$D(A)\hookrightarrow \cal{H}$ is compact. By \cite[Chapter II, Proposition~4.25 p.~117]{EngelNagel} then i) follows.

To prove ii) we remark that the inclusion $\cal{N}\subseteq\text{Ker } A$ trivially follows by Lemma~\ref{lemma 5.1} and \eqref{ad}.
To prove the reverse inclusion we take $U=(u,v,w,z)\in \text{Ker } A$. By \eqref{ad} and \eqref{5.5} we then get that $w=0$, that $z=0$, that
$u\in H^2(\Omega)$ solves the homogeneous Neumann problem
$\Delta u=0$  in $\Omega$, $\partial_\nu u=0$ on $\Gamma$, so $u=c_1\in\C$ in $\Omega$ since $\Omega$ is connected, and finally that $v\in H^2(\Gamma_1)$ solves the elliptic equation
\begin{equation}\label{5.7}
-\DivGamma (\sigma \nabla_\Gamma v)+\kappa v=0 \quad\text{on $\Gamma_1$.}
\end{equation}
Multiplying \eqref{5.7} by $v$ and using \eqref{ridinaform} we then get
$\int_{\Gamma_1}\sigma |\nabla_\Gamma v|_\Gamma^2+\kappa |v|^2=0$. Then $\nabla_\Gamma v=0$ so, being $\Gamma_1$ connected, $v=c_2\in\C$ and, when $\kappa\not\equiv 0$, $c_2=0$.
Hence $U\in \cal{N}$, proving ii).

To prove iii) we remark that, using \eqref{ad}, \eqref{5.1}, \eqref{lap}, \eqref{5.5} and \eqref{ridinaform}, for  all $U=(u,v,w,z)\in D(A)$ we have
\begin{equation*}
		\begin{split}
			[AU,U]_{\mathcal{H}}=&2i\Ima\left\{\rho_0\int_{\Omega}\nabla u\nabla\overline{w}
+\rho_0\int_{\Gamma_1}w\overline{z}+\int_{\Gamma_1}\sigma (\nabla_\Gamma v,\nabla_\Gamma z)_{\Gamma}
+\int_{\Gamma_1}\kappa v\overline{z}
\right\}\\
&+\int_{\Gamma_1}\delta\lvert z\rvert^2\,+\dfrac{\rho_0}{c^2}\int_{\Omega}d\lvert w\rvert^2,
		\end{split}
	\end{equation*}
so \eqref{5.4} immediately follows.

To prove iv) we take $U=(u,v,w,z)\in D(A)$ and $\lambda\in\R$ such that $AU=i\lambda U$. By part ii) we have to prove that
if $d\not\equiv0$ and $\lambda\not=0$ then $U=0$.
By \eqref{5.5} then  $u\in H^2(\Omega)$, $v\in H^2(\Gamma_1)$, $w\in H^1(\Omega)$ and $z\in H^1(\Gamma_1)$ satisfy the system
\begin{equation}\label{5.8}
\begin{cases}
-w =i\lambda u \qquad &\text{in
$\Omega$,}\\
- z=i\lambda v\qquad
&\text{on
$\Gamma_1$,}\\
-c^2\Delta u+dw=i\lambda w \qquad &\text{in
$\Omega$,}\\
- \DivGamma (\sigma \nabla_\Gamma v)+\kappa v+\delta z+\rho_0 w =i\lambda \mu z\qquad
&\text{on
$\Gamma_1$,}\\
\partial_\nu u=z\qquad
&\text{on
$\Gamma_1$,}\\
\partial_\nu u=0 &\text{on $\Gamma_0$,}
\end{cases}
\end{equation}
so in particular $u\in H^2(\Omega)$ solves the elliptic equation
\begin{equation}\label{5.9}
 -c^2\Delta u=(\lambda^2+i\lambda d)u\qquad\text{in $\Omega$.}
\end{equation}
Since $AU=i\lambda U$, trivially $\Real [AU,U]_\cal{H}=0$, so by \eqref{5.4} we have $d|w|^2=0$ in $\Omega$.
Since $r\ge 2$, by assumption (A3) we have $d\in W^{1,\infty}(\Omega)$, so we can fix its representative $d\in C(\Omega)$.
Since $d\not\equiv 0$, $\widetilde{\Omega}:=\{x\in\Omega: d(x)>0\}$ is open and non empty.
Moreover $w=0$ in $\widetilde{\Omega}$. By \eqref{5.8}$_1$, as $\lambda\not=0$,  we then get $u=0$ in $\widetilde{\Omega}$.
Since $u$ solves \eqref{5.9} we can then apply the Unique Continuation Theorem \cite[Theorem, p. 235]{Aronszajn} to get that $u=0$ in $\Omega$.
Then, by \eqref{5.8}$_1$ and \eqref{5.8}$_5$ , $w=0$ in $\Omega$ and $z=0$ on $\Gamma_1$. Since $\lambda\not=0$ by \eqref{5.8}$_2$ we get $v=0$ on $\Gamma_1$,
so $U=0$, completing the proof.
\end{proof}
Lemma~\ref{lemma5.2}--ii) shows that to prove asymptotic stability for the abstract problem \eqref{vectpb} means
proving that for all $U_0\in\cal{H}$ there is $N_0\in\cal{N}$ such that $T(t)[U_0]\to N_0$ as $t\to\infty$.

Lemma~\ref{lemma5.2} also suggests the following general strategy in order to prove such a stability result.
We shall look for a closed subspace $\cal{M}$ of $\cal{H}$ with the following properties:
\renewcommand{\labelenumi}{{\alph{enumi})}}
\begin{enumerate}
\item $\cal{M}$ is invariant under the flow of $\{T(t),t\in\R\}$;
\item $\cal{M}$ complements $\cal{N}$, that is $\cal{H}=\cal{N}\oplus\cal{M}$;
\item the restrictions of $\n\cdot\n_\cal{H}$ and $\|\cdot\|_\cal{H}$ to $\cal{M}$ are two equivalent norms  on it.
\end{enumerate}
The natural starting point to find $\cal{M}$ is the property a). To satisfy the property b) the subspace $\cal{M}$ has to depend
on the alternative $\kappa\equiv 0$ versus $\kappa\not\equiv 0$, since $\cal{N}$ depends on it by Lemma~\ref{lemma 5.1}.
Hence we shall look for a  functional invariant under the flow and, when $\kappa\equiv 0$, for an additional one.

The first of them is nothing but the trivial generalization to the case $d\ge 0$ of the functional given in \cite[Lemma~6.2.1]{mugnvit},
although it plays a different role in the quoted paper. The second one is directly suggested by equation \eqref{1.1}$_2$.
We set $L_1,L_2\in\cal{H}'$ by
\begin{equation}\label{5.10}
  L_1(u,v,w,z)=\int_\Omega w+\int_\Omega du-c^2\int_{\Gamma_1}v,\qquad
  L_2(u,v,w,z)=\int_{\Gamma_1} \mu z+\delta v+\rho_0 u,
\end{equation}
and
\begin{equation}\label{5.11}
  \cal{M}=
  \begin{cases}
  \text{Ker }L_1, \quad & \text{if $\kappa\not\equiv 0$,}
  \\
  \text{Ker }L_1\cap \text{Ker }L_2, \quad & \text{if $\kappa\equiv 0$.}
\end{cases}
\end{equation}
The following result holds.
\begin{lemma}\label{Lemma 5.3}
Let (A0--4) hold. Then for all $U_0\in\cal{H}$
\begin{equation}\label{5.12}
  L_1(T(t)[U_0])=L_1[U_0]\qquad\text{for all $t\in\R$,}
\end{equation}
and, when $\kappa\equiv 0$,
\begin{equation}\label{5.13}
  L_2(T(t)[U_0])=L_2[U_0]\qquad\text{for all $t\in\R$.}
\end{equation}
Consequently $\cal{M}$ is invariant under the flow of $\{T(t), t\in\R\}$.
\end{lemma}
\begin{proof}
When $U_0\in D(A)$, by \eqref{ad}, \eqref{vectpb}  and \eqref{5.5} we have
\begin{align*}
&\dfrac d{dt}\left( \int_\Omega w+\int_\Omega du-c^2\int_{\Gamma_1} v\right)=
\int_\Omega w_t+\int_\Omega dw-c^2\int_{\Gamma_1} z=c^2\left(\int_\Omega \Delta u-\int_\Gamma \partial_\nu u\right)
\\ \intertext{and}
&\dfrac d{dt} \int_{\Gamma_1} \mu z+\delta v+\rho_0 u=
\int_{\Gamma_1} \mu z_t+\delta v_t+\rho_0u_t=\int_{\Gamma_1} \DivGamma (\sigma \nabla_\Gamma v)-\kappa v.
\end{align*}
Hence $\frac d{dt} L_1(T(t)[U_0])=0$ by the Divergence Theorem in $H^2(\Omega)$, and when $\kappa\equiv 0$,
$\frac d{dt} L_2(T(t)[U_0])=0$ by \eqref{ridinaform}. Since $L_i(T(\cdot)[U_0])\in C^1(\R)$ for $i=1,2$, by integrating on time we get \eqref{5.12} and \eqref{5.13} when $\kappa\equiv 0$.  We then extend them to data in $\cal{H}$ by density.
\end{proof}
The following result shows that $\cal{M}$ also satisfy the property b) when $d\not\equiv 0$. When $d\equiv 0$ as in \cite{mugnvit} one has $L_1(1,0,0,0)=0$,
so  when also $\kappa\not\equiv 0$ one has $\cal{M}\cap \cal{N}\not=\{0\}$ and such an approach is impossible. To get the precise asymptotic in Theorem~\ref{stabt}
we also make explicit the projectors $\Pi_\cal{N}$ and  $\Pi_\cal{M}$.
\begin{lemma}\label{Lemma 5.4}
Let (A0--4) hold and $d\not\equiv 0$. Then the following splitting
\begin{equation}\label{5.14}
  \cal{H}=\cal{N}\oplus\cal{M}
\end{equation}
holds. The associated projectors $\Pi_\cal{N}\in\cal{L}(\cal{H},\cal{N})$ and  $\Pi_\cal{M}\in\cal{L}(\cal{H},\cal{M})$ are given,
for all $U=(u,v,w,z)\in \cal{H}$, as follows.
\renewcommand{\labelenumi}{{\roman{enumi})}}
\begin{enumerate}
\item When $\kappa\not\equiv 0$, using $\alpha$ given in \eqref{alpha},
\begin{equation}\label{5.17}
\Pi_\cal{N} U=(\alpha U, 0,0,0),\qquad \Pi_\cal{M} U=(u-\alpha U, v,w,z).
\end{equation}
\item When $\kappa\equiv 0$, using  $\beta$ and $\gamma$  given by \eqref{beta},
\begin{equation}\label{5.18}
\Pi_\cal{N} U=(\beta U, \gamma U,0,0),\qquad \Pi_\cal{M} U=(u-\beta U, v-\gamma U,w,z).
\end{equation}
\end{enumerate}
\end{lemma}
\begin{proof}
We first consider the case $\kappa\not\equiv 0$, in which by \eqref{5.3} and \eqref{5.11}
\begin{equation}\label{5.19}
 \cal{N}=\C V_1,\quad\text{and}
 \quad \cal{M}=\text{Ker }L_1,\quad\text{where}\quad V_1:=(1,0,0,0).
\end{equation}
Since $d\ge 0$ and $d\not\equiv 0$ one has $L_1 V_1=\int_\Omega d>0$, and consequently the only solution $x\in\C$ of the equation $L_1(xV_1)=0$ is the trivial one, so proving that
\begin{equation}\label{5.20}
\cal{N}\cap\cal{M}=\{0\}.
\end{equation}
Moreover, for all $U\in\cal{H}$, setting
$U'=\frac {L_1U}{L_1V_1}V_1$, and $U''=U-U'$,
one trivially has $U'\in\cal{N}$ and $L_1U''=L_1U-L_1U=0$, so $U''\in\cal{M}$. Hence $\cal{H}=\cal{N}+\cal{M}$, which combined with \eqref{5.20} gives \eqref{5.14}. Setting $\Pi_\cal{N}U=U'$ and $\Pi_\cal{M}U=U''$ (see \cite[Chapter~I,~\S~4,~p. 20]{Kato}), the proof in the case $\kappa\not\equiv 0$ is complete.

We now consider the case $\kappa\equiv 0$, in which, $V_1$ being still defined by \eqref{5.19},
\begin{equation}\label{5.22}
 \cal{N}=\text{span }\{V_1,V_2\}\quad\text{and}
 \quad \cal{M}=\text{Ker }L_1\cap\text{Ker }L_2,\quad\text{where}\quad V_2:=(0,1,0,0).
\end{equation}
Since
$$L_1V_1=\int_\Omega d,\quad L_1V_2=-c^2 \cal{H}^{N-1}(\Gamma_1),\quad L_2V_1=\rho_0\cal{H}^{N-1}(\Gamma_1), \quad L_2V_2=\int_{\Gamma_1}\delta,$$
the matrix
$$C:=\left(
    \begin{array}{cc}
      L_1V_1 & L_1V_2 \\
      L_2V_1 & L_2V_2 \\
    \end{array}
  \right)$$
  has determinant
$$\text{det }C=\int_\Omega d\int_{\Gamma_1}\delta +\rho_0c^2 [\cal{H}^{N-1}(\Gamma_1)]^2>0.$$
Consequently the linear homogeneous system $L_1(xV_1+yV_2)=L_2(xV_1+yV_2)=0$ has only the trivial solution $x=y=0$, so by \eqref{5.22} also in this case \eqref{5.20} holds true. Moreover, for the same reason, for all $U\in\cal{H}$, its nonhomogeneous version
\begin{equation}\label{5.27}
 L_1(xV_1+yV_2)=L_1U, \qquad L_2(xV_1+yV_2)=L_2U
\end{equation}
 has a unique solution $(x,y)\in\C^2$, which by Cramer's formula is $x=\beta U$ and $y=\gamma U$, where $\beta$ and $\gamma$ are given by \eqref{beta}. We then set $U'=xV_1+yV_2$ and $U''=U-U'$. Trivially $U'\in\cal{N}$ and, by \eqref{5.27}, $L_1U''=L_2U''=0$, so $U''\in\cal{M}$. Hence $\cal{H}=\cal{N}+\cal{M}$, which combined with \eqref{5.20} gives \eqref{5.14} also in the case $\kappa\equiv 0$. We then complete the proof as in the previous case.
\end{proof}
To prove that $\cal{M}$ also enjoys the property c) we recall the following probably well--known result, a proof of which was given in \cite[Lemma~6.1.2]{mugnvit}.
\begin{lemma}\label{Lemma 5.5}Let $[\cdot,\cdot]$ be a continuous pseudo--inner product on a real or complex Hilbert space $(H, (\cdot,\cdot))$, and let $N=\{u\in H: [u,u]=0\}$ be finite dimensional. Then $[\cdot,\cdot]$ is coercive on any closed subspace $H_1$ of $H$ such that $H_1\cap N=\{0\}$ if and only if
it is coercive on $N^\bot$.
\end{lemma}
By using Lemmas~\ref{Lemma 5.4} and Lemma~\ref{Lemma 5.5} we then get
\begin{lemma}\label{Lemma 5.6}
If (A0--4) hold and $d\not\equiv 0$ the restriction of $\n\cdot\n_{\cal{H}}$ to $\cal{M}$ is a norm equivalent to $\|\cdot\|_{\cal{H}}$ on it.
\end{lemma}
\begin{proof}
By Lemmas~\ref{Lemma 5.4} and Lemma~\ref{Lemma 5.5} we only have to show that $\n\cdot\n_{\cal{H}}$ is coercive on $\cal{N}^\bot$. By \eqref{innerp} and \eqref{5.3}, since $\rho(x)=\rho_0>0$, one immediately gets
\begin{equation}\label{5.32}
\cal{N}^\bot=
\begin{cases}
 H^1_c(\Omega)\times H^1(\Gamma_1)\times L^2(\Omega)\times L^2(\Gamma_1)&\quad\text{when $\kappa\not\equiv 0$},\\
  H^1_c(\Omega)\times H^1_c(\Gamma_1)\times L^2(\Omega)\times L^2(\Gamma_1)&\quad\text{when $\kappa\equiv 0$,}
\end{cases}
\end{equation}
where we denote
$$H_c(\Omega)=\left\{u\in H^1(\Omega):\,\int_\Omega u=0\right\},\quad\text{and}\quad
H_c(\Gamma_1)=\left\{v\in H^1(\Gamma_1):\,\int_{\Gamma_1} v=0\right\}.$$
By respectively setting on $H^1_c(\Omega)$, $H^1(\Gamma_1)$, $L^2(\Omega)$ and $L^2(\Gamma_1)$ the pseudo--inner products
\begin{alignat*}2
&[u_1,u_2]_{H^1_c(\Omega)}=\rho_0\int_\Omega \nabla u_1\nabla\overline{u_2},\qquad
&& [v_1,v_2]_{H^1(\Gamma_1)}=\int_{\Gamma_1}\sigma(\nabla_\Gamma v_1,\nabla_\Gamma v_2)_{\Gamma}+\int_{\Gamma_1}\kappa v_1\overline{v_2},\\
&[w_1,w_2]_{L^2(\Omega)}=\frac{\rho_0}{c^2}\int_\Omega w_1\overline{w_2},\qquad
&&[z_1,z_2]_{L^2(\Gamma_1)}=\int_{\Gamma_1}\mu z_1\overline{z_2},
\end{alignat*}
we make the following remarks. At first $[\cdot,\cdot]_{H^1_c(\Omega)}$ is coercive on $H^1_c(\Omega)$ by a Poincar\`e type inequality. See for example
\cite[Chapter~13,~Theorem~13.2.7,~p.423]{LeoniSobolev2}. Next  $[\cdot,\cdot]_{H^1(\Gamma_1)}$ is coercive on $H^1(\Gamma_1)$ when $\kappa\not\equiv 0$ and on
$H^1_c(\Gamma_1)$ when $\kappa\equiv 0$ by Poincar\`e type inequalities. See \cite[Lemma~6.1.4]{mugnvit}. Finally  $[\cdot,\cdot]_{L^2(\Omega)}$
and $[\cdot,\cdot]_{L^2(\Gamma_1)}$ are coercive since $\rho_0,c>0$ and $\mu\ge \mu_0>0$.
Hence, by \eqref{5.1} and \eqref{5.32}, $\n\cdot\n_{\cal{H}}$ is coercive on $\cal{N}^\bot$, concluding the proof.
\end{proof}
We can finally give
\begin{proof}[Proof of Theorem~\ref{stabt}]
By Lemma~\ref{Lemma 5.3} and standard semigroup theory, see \cite[Chapter~I,~\S5.11,~p.43]{EngelNagel}, the restrictions of the operators $T(t)$ to $\cal{M}$
for $t\in\R$ constitute the strongly continuous subspace group $\{T_{\cal M}(t),t\in\R\}$. By Lemma~\ref{Lemma 5.6} we can equivalently equip $\cal{M}$
with the restriction of $[\cdot,\cdot]_\cal{H}$. By \cite[Chapter~II,~\S2.3,~Corollary,~p.61]{EngelNagel} and Theorem~\ref{almostwell}--i) the generator of
$\{T_{\cal M}(t),t\in\R\}$ is the operator $-A_\cal{M}$, where
$$A_\cal{M}:D(A_\cal{M})\subset\cal{M}\to\cal{M},\quad D(A_\cal{M})=D(A)\cap\cal{M},\quad A_\cal{M}=A_{|\cal{M}}.$$
Moreover, by Lemma~\ref{lemma5.2}--ii), the semigroup $\{T_{\cal M}(t),t\ge 0\}$ is contractive and then bounded. By \cite[Chapter~II,~Theorem~1.10,~p.55]{EngelNagel}
the resolvent of $-A_\cal{M}$ is nothing but the restriction to $\cal{M}$ of the resolvent of $-A$ and then, by Lemma~\ref{lemma5.2}--i),
$-A_\cal{M}$ has compact resolvent. Consequently, by \cite[Chapter~IV,~Corollary~1.19,~p.248]{EngelNagel}, the spectrum of $A_\cal{M}$ reduces to its point spectrum
$\sigma_p(A_\cal{M})$ and, see \cite[Chapter~III,~p.187]{Kato}, is countable.
Moreover, by Lemma~\ref{Lemma 5.4}, $\cal{N}\cap\cal{M}=\{0\}$. Hence, by Lemma~\ref{lemma5.2}--ii),
$\text{Ker }A_\cal{M}=\text{Ker }A\cap\cal{M}=\cal{N}\cap\cal{M}=\{0\}$. By Lemma~\ref{lemma5.2}--iv) we then have
$\sigma_p(A_\cal{M})\cap i\R=\emptyset$.

We can then apply the consequence of the  Arendt, Batty, Lyubich and V\~{u} Theorem for reflexive space,
see \cite[Chapter V, Corollary~2.22 p.~327]{EngelNagel}), to get that the semigroup $\{T_{\cal M}(t),t\ge 0\}$ is strongly stable, that is
\begin{equation}\label{5.35}
 T_\cal{M}(t)[U_0]\to 0,\quad\text{in $\cal{M}$, \quad as $t\to\infty$, \quad for all $U_0\in\cal{M}$}.
\end{equation}
Since $\cal{N}=\text{Ker }A$, the subspace group of $\{T(t),t\in\R\}$ on $\cal{N}$ reduces to the identity. Consequently, by Lemma~\ref{Lemma 5.4},
we have
\begin{equation}\label{5.36}
  T(t)=T_\cal{M}(t)\cdot \Pi_\cal{M}+\Pi_\cal{N},\qquad\text{for all $t\in\R$.}
\end{equation}
Combining \eqref{5.35} and \eqref{5.36} we get
$$T(t)[U_0]\to \Pi_\cal{N}U_0,\quad\text{in $\cal{H}$, \quad as $t\to\infty$, \quad for all $U_0\in\cal{H}$}.$$
By combining it with \eqref{5.17} and \eqref{5.18} we get \eqref{S1} and \eqref{S2}.  Combining them with \eqref{energyfunctional} we conclude the proof.
\end{proof}
\bigskip
{\bf Statements and Declarations.} The authors declare that they have  no
conflict of interest.
%\bibliographystyle{amsplain}
%\bibliography{biblio2}
\def\cprime{$'$}
\providecommand{\bysame}{\leavevmode\hbox to3em{\hrulefill}\thinspace}
\providecommand{\MR}{\relax\ifhmode\unskip\space\fi MR }
\providecommand{\MRhref}[2]{%
  \href{http://www.ams.org/mathscinet-getitem?mr=#1}{#2}
}
\providecommand{\href}[2]{#2}

\end{document}